\documentclass{amsart}
\usepackage{amscd,amssymb}
\usepackage[dvips]{graphics}
\input diagrams
\diagramstyle[PostScript=dvips]

\newcommand{\includefig}[1]{\raisebox{-3ex}{\resizebox{!}{7ex}{\includegraphics{#1.eps}}}}
\newcommand{\tcongi}{\includefig{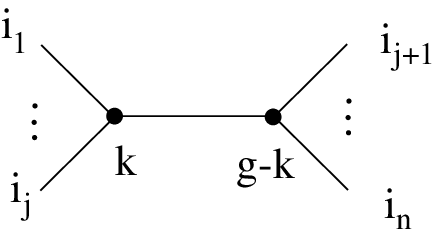}}
\newcommand{\tcongm}{\includefig{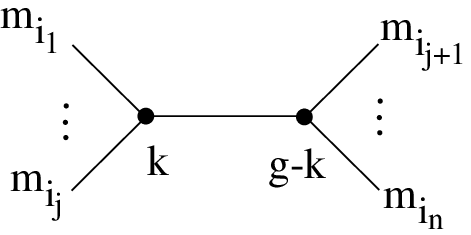}}
\newcommand{\ocongi}{\includefig{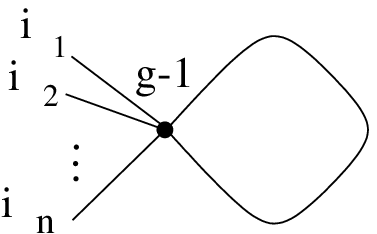}}

\newcommand{\scheme}[1]{\mathcal {#1}}

\newcommand{\mbar}{\M}

\newcommand{\mgn}{\scheme{M}_{g,n}}
\newcommand{\mgnrm} {\mgn^{1/r,\bm}}

\newcommand{\mgnbar}{\overline{\scheme{M}}_{g,n}}
\newcommand{\mgnrmbar} {\mgnbar^{1/r,\bm}}
\newcommand{\mgnrmvbbar} {\mgnbar^{1/r,\bm}(V,\beta)}
\newcommand{\mgnvbbar} {\mgnbar(V,\beta)}
\newcommand{\mgnrvbbar} {\mgnbar^{1/r}(V,\beta)}
\newcommand{\mgammabar}{\M_{\Gamma}}
\newcommand{\mgammavbbar} {\mgammabar(V,\beta)}
\newcommand{\mgammarbar} {\mgammabar^{1/r}}
\newcommand{\mgammarvbbar} {\mgammarbar(V,\beta)}
\newcommand{\mgnrbar} {\mgnbar^{1/r}}

\newcommand{\cgn}{\scheme{C}_{g,n}}

\newcommand{\cgnrm}{\cgn^{1/r,\bm}}

\newcommand{\pic}{\text{Pic\,}}        

\newcommand{\sym}{\mathbf{Sym}}
\newcommand{\tensor}{\otimes}
\newcommand{\cross}{\times}
\newcommand{\irightarrow}{\rTo^{\sim}}

\newcommand{\ck}{{\mathcal K}}
\newcommand{\cl}{{\mathcal L}}
\newcommand{\ce}{{\mathcal E}}
\newcommand{\cf}{{\mathcal F}}
\newcommand{\co}{{\mathcal O}}
\newcommand{\cc}{{\mathcal C}}

\newcommand{\ch}{{\mathcal H}}

\newcommand{\cR}{{\mathcal R}}

\newcommand{\proj}{\operatorname{Proj}}
\newcommand{\spec}{\operatorname{Spec}}
\newcommand{\maxid}{\mathfrak{m}}
\newcommand{\fp}{\mathfrak{p}}
\newcommand{\fq}{\mathfrak{q}}

\newcommand{\be}{\mathbf{e}}
\newcommand{\bev}{\mathbf{ev}}
\newcommand{\bevt}{\tilde{\mathbf{ev}}}
\newcommand{\bdelta}{\boldsymbol{\delta}}
\newcommand{\bgamma}{\boldsymbol{\gamma}}
\newcommand{\bm}{\mathbf{m}}          

\newcommand{\bmt}{\widetilde{\mathbf{m}}}         

\newcommand{\bone}{\mathbf{1}}

\newcommand{\bt}{\mathbf{t}}
\newcommand{\btau}{\boldsymbol{\tau}}
\newcommand{\btaut}{\widetilde{\btau}}
\newcommand{\btT}{\tilde{\boldsymbol{t}}}

\newcommand{\bu}{\mathbf{u}}
\newcommand{\bx}{\mathbf{x}}       
\newcommand{\by}{\mathbf{y}}       
\newcommand{\cft}{\mbox{CohFT}}              
\newcommand{\cfts}{\mbox{CohFTs}}           
 
\newcommand{\chr}{{\ch^{(r)}}}

\newcommand{\chVr}{{\ch^{(V,r)}}}      
\newcommand{\cor}[1]{\langle\,{#1}\,\rangle}

\newcommand{\cv}{c^{1/r}}           
\newcommand{\cvt}{\tilde{c}^{1/r}}        

\newcommand{\cvtim}{\tilde{c}^{1/r,\bm,(i)}}
\newcommand{\ee}{\varepsilon}        
\newcommand{\ev}{ev}         
\newcommand{\evt}{\tilde{ev}}
\newcommand{\etar}{{\eta^{(r)}}}    
\newcommand{\etaVr}{{\eta^{(V,r)}}} 
\newcommand{\etaV}{\eta}

\newcommand{\KdV}{\mathrm{KdV}}
\newcommand{\Lambdar}{\Lambda^{(r)}} 
 
\newcommand{\LambdaV}{\Lambda^{(V)}}    
\newcommand{\LambdaVr}{\Lambda^{(V,r)}} 
\newcommand{\M}{\overline{\MM}}     
\newcommand{\MM}{\scheme{M}}         
\newcommand{\mt}{{\tilde{m}}}        
\newcommand{\nc}{{\mathbb{C}}}        
\newcommand{\nq}{{\mathbb{Q}}}        
        
\newcommand{\nz}{{\mathbb{Z}}}        
\newcommand{\PhiV}{\Phi^{(V)}}

\newcommand{\PhiVr}{\Phi^{(V,r)}}
\newcommand{\PhitVr}{\widetilde{\Phi}^{(V,r)}}
\newcommand{\pt}{\tilde{p}}
\newcommand{\pti}{\pt_{(i)}}
\newcommand{\st}{\mathrm{st}}         
\newcommand{\taut}{\tilde{\tau}}
\newcommand{\tT}{\tilde{t}}

\newcommand{\cpone}{\mathbb{P}^1}

\newtheorem{thm}{Theorem}[section]
\newtheorem{lm}[thm]{Lemma}
\newtheorem{prop}[thm]{Proposition}
\newtheorem{crl}[thm]{Corollary}

\theoremstyle{definition}

\newtheorem{rem}[thm]{Remark}

\newtheorem{df}[thm]{Definition}
\newtheorem{ex}[thm]{Example}
\newtheorem{df-pr}[thm]{Definition-Proposition}

\theoremstyle{remark}
 
\newenvironment{ack}{\begin{sloppypar}\emph{Acknowledgments}.}{\end{sloppypar}}

\begin{document}
\addtocounter{section}{-1}

\title[Stable Spin Maps]
{Stable Spin Maps, Gromov-Witten Invariants, and Quantum
Cohomology}

\subjclass{Primary: 14N35, 53D45. Secondary: 14H10}

\author
[T. J. Jarvis]{Tyler J. Jarvis}
\address
{Department of Mathematics, Brigham Young University, Provo, UT
84602, USA} \email{jarvis@math.byu.edu}
\thanks{Research of the first author was partially supported by NSA
grant number MDA904-99-1-0039}

\author
[T. Kimura]{Takashi Kimura}
\address
{Department of Mathematics, 111 Cummington Street, Boston
University, Boston, MA 02215, USA} \email{kimura@math.bu.edu}
\thanks{Research of the second author was partially supported by NSF grant
  number DMS-9803427}

\author
[A. Vaintrob]{Arkady Vaintrob}
\address
{Department of Mathematics, University of Oregon, Eugene, OR
974003, USA} \email{vaintrob@math.uoregon.edu}

\date{}

\begin{abstract}

We introduce the stack $\M_{g,n}^{1/r}(V)$ of $r$-spin maps. These
are stable maps into a variety $V$ from $n$-pointed algebraic
curves of genus $g$, with the additional data of an $r$-spin
structure on the curve.

We prove that $\M_{g,n}^{1/r}(V)$ is a Deligne-Mumford stack, and
we define analogs of the Gromov-Witten classes associated to these
spaces. We show that these classes yield a cohomological field 
theory (\cft) that is the tensor product of the \cft\ associated 
to the usual Gromov-Witten invariants of $V$ and the $r$-spin 
\cft. When $r=2$, our construction gives the usual Gromov-Witten 
invariants of $V$.

Restricting to genus zero, we obtain the notion of an $r$-spin
quantum cohomology of $V$, whose Frobenius structure is isomorphic
to the tensor product of the Frobenius manifolds corresponding to 
the quantum cohomology of $V$ and  the $r$-th Gelfand-Dickey 
hierarchy (or, equivalently, the $A_{r-1}$ singularity). We also 
prove a generalization of the descent property which, in 
particular, explains the appearance of the $\psi$ classes in the 
definition of gravitational descendants. Finally, we compute the 
small phase space potential function when $r=3$ and $V=\cpone$. 
\end{abstract}

\maketitle

\section{Introduction}

In this paper, we present a generalization of the theory of
quantum cohomology and Gromov-Witten invariants arising from
algebraic curves with higher spin structures. Recall that ordinary 
Gromov-Witten invariants of a projective variety  $V$ are 
constructed by means of the moduli spaces $\M_{g,n}(V)$ of stable 
maps to $V$. The space $\M_{g,n}(V)$ is a Deligne-Mumford stack 
compactifying the space of holomorphic maps to $V$ from Riemann 
surfaces of genus $g$  with $n$ marked points. In particular, the 
moduli space of stable maps to a point  coincides with the moduli 
of stable curves $\M_{g,n}$. 

Although the space $\M_{g,n}(V)$ is not generally smooth, it does
possess a homology class $[\M_{g,n}(V)]^\mathrm{virt}$ (called the 
virtual fundamental class) which plays the role of the fundamental 
class in intersection theory. Furthermore, the  evaluation maps 
$\M_{g,n}(V)\to V$ obtained by evaluating the  map on each marked 
point on the curve allow one to pull back cohomology classes from 
$V$. Together, these classes give rise to the collection 
$\Lambda^V\,\:=\,\{\,\Lambda^V_{g,n}\,\}$ of Gromov-Witten classes 
$$ 
 \Lambda^V_{g,n}\in H^\bullet(\M_{g,n})\,\otimes\,T^n H^\bullet(V)^*,
$$ which behave nicely when restricted to the boundary strata of
$\M_{g,n}$. This allows one to define a collection of multilinear
operations on the space $H^\bullet(V)$, parametrized by elements
of $H_\bullet(\M_{g,n})$. These operations satisfy the axioms of a
cohomological field theory (\cft) in the sense of
Kontsevich-Manin~\cite{KM1}. In particular, their restriction to
stable maps of genus zero endows $H^\bullet(V)$ with the structure
of a (formal) Frobenius manifold~\cite{Du,Hi,Ma}, called the
quantum cohomology of $V$, whose multiplication is a deformation 
of the usual cup product in $H^\bullet(V)$.

The diagonal map $\M_{g,n}\,\to\,\M_{g,n}\,\times\,\M_{g,n}$
induces  a natural tensor product operation on the set of \cfts.
Behrend~\cite{B} proved that the tensor product of the \cfts\
defined by the Gromov-Witten classes of two smooth projective
varieties $V$ and $V'$ is isomorphic to the \cft\ defined by the
classes $\Lambda^{V\,\times\,V'}$. This is nontrivial because
$\M_{g,n}(V\,\times\,V')$ is not isomorphic to the fiber product
$\M_{g,n}(V)\,\times_{\M_{g,n}}\,\M_{g,n}(V')$. Restricting to 
genus zero, we can regard this result as a K\"{u}nneth formula for 
quantum cohomology. 
\medskip

In~\cite{JKV}, we   introduced a new class of \cfts, one for each
integer $r\geq 2$, based on the moduli space \ \ $\M_{g,n}^{1/r}=
\displaystyle \coprod_{\bm} \M_{g,n}^{1/r,\bm}$ \ \ of higher spin
curves, constructed in~\cite{J}. Recall that for 
$\bm\,=\,(m_1,\ldots,m_n)$, with $m_i \in \nz$,  the moduli space  
$\M_{g,n}^{1/r,\bm}$ is a compactification of the space of Riemann 
surfaces of genus $g$ with $n$ marked points $p_1,\ldots,p_n$ and 
an $r$-th root of the twisted canonical line bundle $\omega 
\otimes \co(-\sum_i m_i p_i)$. 

We proved in~\cite{JKV} that the existence of a special cohomology
class $\cv$ (called the spin virtual class) in
$H^\bullet(\M_{g,n}^{1/r})$, satisfying certain axioms similar to
the Behrend-Manin axioms~\cite{BMa} for the virtual fundamental
class, gives rise to a \cft\ of rank $r-1$.

This $r$-spin \cft\ is related to the work of Witten~\cite{W}, who
conjectured that its large phase space potential (the generating
function of certain intersection numbers on $\M_{g,n}^{1/r,\bm}$)
is a $\tau$ function of the \mbox{$r$-th} Gelfand-Dickey (or
$\KdV_r$) hierarchy. When $r=2$, this conjecture reduces to an
earlier conjecture of Witten's on the intersection numbers of
$\M_{g,n}$, which was proved by Kontsevich~\cite{Ko}.
In~\cite{JKV},  following Witten's ideas, we constructed the spin
virtual class and proved the conjecture in the cases $g=0$ (for
all $r$) and $r=2$ (for all $g$).

In this paper, we introduce and describe moduli spaces that give
an intersection-theoretic realization of the tensor product of the
Gromov-Witten \cft\ and the $r$-spin \cft. We construct
$\M_{g,n}^{1/r}(V)$, the stack of stable $r$-spin maps into a
projective variety $V$---objects which combine both the data of a
stable map and an $r$-spin structure.  We prove that
$\M_{g,n}^{1/r}(V)$ is a Deligne-Mumford stack and a ramified
cover of $\M_{g,n}(V)$. Similar to the case of ordinary stable 
maps, stable $r$-spin maps to a point are just stable $r$-spin 
curves.

We introduce the class $\cvt\in H^\bullet(\M_{g,n}^{1/r}(V))$,
which is an analog of a virtual class $\cv$ on $\M_{g,n}^{1/r}$.
The class $\cvt$ is defined by pulling back $\cv$ from
$\M_{g,n}^{1/r}$ when $2g-2+n>0$ and is defined by a direct
construction  for other values of $g$ and $n$. Using the class
$\cvt$ and the virtual fundamental class of $\M_{g,n}(V)$, we
define the $r$-spin analogs of the Gromov-Witten classes $$
\LambdaVr_{g,n}\,\in\,H^\bullet(\M_{g,n}^{1/r}(V))\,\otimes\, T^n
{\chVr}^*, $$ where $\chVr\,=\,H^\bullet(V)\otimes{\chr}$ and
$\chr$ is the state space of the $r$-spin \cft. The collection of
classes $\LambdaVr = \{\,\LambdaVr_{g,n}\,\}$ gives rise to a
\cft\ with the state space $\chVr$, which is isomorphic to the
tensor product of the Gromov-Witten \cft\ and the $r$-spin \cft.
This fact is not trivial because the space $\M_{g,n}^{1/r}(V)$ is
not isomorphic to the fiber product
$\M_{g,n}(V)\,\times_{\M_{g,n}}\,\M_{g,n}^{1/r}$.
\medskip

Our result should be viewed as an analog of the relationship
between stable maps to $V\,\times\,V'$ and the tensor product of
the \cfts\ corresponding to each factor (see~\cite{KM1,KMK,B}).
Restricting to genus zero,   we obtain the tensor product of the
quantum cohomology of $V$ with the Frobenius manifold associated
to $\KdV_r$ (or equivalently, to the $A_{r-1}$ singularity).
Tensor products of several spin \cfts\ were studied in our earlier
paper~\cite{JKV2}. They correspond to the moduli spaces of curves
with multiple spin structures. It would be very interesting to
find an enumerative interpretation of the corresponding invariants
similar to the interpretation of the ordinary Gromov-Witten
invariants.

Finally, it is worth observing that these results have a physical
interpretation as well.              Our spin Gromov-Witten
invariants may be regarded as the correlators in a theory of
topological gravity coupled to topological matter, where the
matter sector of the theory is the topological sigma model with
target space $V$ coupled with a certain type of gauged
$SU(2)_{r-2}/U(1)$ Wess-Zumino-Witten model. Our moduli spaces
$\M_{g,n}^{1/r}(V)$ provide an intersection-theoretic realization
of this theory.
\medskip

The structure of the paper is as follows. In the first section, we
recall some standard definitions and properties of stable maps and
Gromov-Witten invariants. In the second section, we construct the
moduli space $\M_{g,n}^{1/r}(V)$, establish some of its
properties, and prove that it is a Deligne-Mumford stack. In the
third section, we introduce decorated graphs associated with the
various moduli spaces and the corresponding cohomology classes.
 In the fourth section, we define the
spin virtual class $\tilde{c}^{1/r}$ on $\M_{g,n}^{1/r}(V)$ and
study its intersection-theoretic properties. In the fifth section,
we define the analogs of the Gromov-Witten classes associated to
$\M_{g,n}^{1/r}(V)$ and show that the resulting \cft\ is the
tensor product of \cfts\ associated to $\M_{g,n}(V)$ and
$\M^{1/r}_{g,n}$. We then define the potential functions of the
theory
 and prove that when $r=2$, this theory reduces to the
usual Gromov-Witten invariants of $V$. Furthermore, we prove that
$\cvt$ satisfies the descent property in genus zero. In the final
section, we compute the genus zero small phase space potential
function associated to $\M_{g,n}^{1/3}(\cpone)$.

\begin{ack}
Parts of the paper were written while T.K.\ was visiting the
Universit\'{e} de Bourgogne  and A.V.\ was visiting Institute des
Hautes \'Etudes Scientifiques. We would like to thank these
institutions for hospitality and support.
\end{ack}

\section{Preliminaries}

In this section we recall the basic definitions related to stable
maps and Gromov-Witten invariants. For details see \cite{KM1,Ma}.

\subsection{Prestable curves}
\

By a {\em curve\/}, in this paper, we mean a reduced, complete,
connected, one-dimensional scheme over $\nc$.  By {\em genus\/} of
a curve $X$, we mean its arithmetic genus  $g=\dim H^1(X,\co_X)$.
An {\em $n$-pointed, prestable curve\/} is a curve $X$ with $n$
distinct marked points $p_1,\ldots,p_n \in X$, such that $X$ has
at most ordinary double points (nodes) as singularities and
$p_1,\ldots,p_n$ are nonsingular.  A prestable curve $X$ is called
{\em stable\/} if, in addition, each irreducible component of $X$
of genus $0$ has at least three distinguished points, and each
component of genus 1 has at least one distinguished point, where a
distinguished point is either a node or a marked point.

 A    {\em family of prestable, $n$-pointed curves\/}
is a flat, proper morphism $ X \rightarrow T $  with $n$ sections
$\fp_1,\ldots,\fp_n : T \to X$, such that each geometric fiber
$(X_t,\fp_1(t),\ldots,\fp_n(t))$ is an $n$-pointed prestable
curve.

By a {\em line bundle\/} we mean an invertible (locally free of
rank one) coherent sheaf. The \emph{canonical sheaf}  of a family
of curves $f:X \rightarrow T$ is the relative dualizing sheaf of
$f$. This sheaf will be denoted by $\omega_{X/T}$, or $\omega_f$.
When $T=\spec \nc$,  we will write $\omega_X$ or  $\omega$ instead
of $\omega_{X/T}$.  Note that for a family of prestable curves,
the canonical sheaf is a line bundle.

A {\em relatively torsion-free sheaf\/} (or just torsion-free
sheaf) on a family of prestable curves $f: X \rightarrow T$ is a
coherent $\co_{X}$-module $\ce$ that is flat over $T$, such that
on each fiber $X_t =X \times_T \spec k(t)$ the restriction $\ce_t$
has no associated primes of height one.  We will only be concerned
with rank-one  torsion-free sheaves.  Such sheaves are sometimes
called {\em admissible}  or {\em sheaves of pure dimension $1$}.
On the open set where $f$ is smooth, a torsion-free sheaf is
locally free.

\subsection{Moduli of stable maps and virtual fundamental class}
\label{sec:fundcl}
\

Let $V$ be   an algebraic variety over $\nc$. A {\em stable map\/} 
to $V$ consists of an $n$-pointed, prestable curve $(X, 
p_1,\ldots,p_n)$ and a morphism $f: X \to V$, such that each 
irreducible component $C$ of $X$ mapped by $f$ to a point is 
stable (i.e.\ $C$ has at least three distinguished points if it is 
of genus zero and at least one distinguished point if it is of 
genus one). A {\em family of stable maps\/} to $V$ is a family of 
prestable, $n$-pointed curves $\pi: X\to T$ with a morphism 
$f:X\to V$, such that the restriction of $f$ to each geometric 
fiber of $\pi$ is stable.

For $\beta\in H_2(V,\nz)$, we say that the stable map $f:X\to V$
has {\em class\/} $\beta$ if the pushforward $f_*([X])$ of the
fundamental class of $X$ is equal to $\beta$. We denote by
$\M_{g,n}(V,\beta)$ the stack of $n$-pointed stable maps of genus
$g$ to $V$ of class $\beta$, and by $\M_{g,n}(V)$ the disjoint
union $\displaystyle\coprod_{\beta\in
H_2(V,\nz)}\M_{g,n}(V,\beta)$. If $V$ is a projective variety,
$\M_{g,n}(V,\beta)$ is a Deligne-Mumford stack. When $V$ is a
point, $\M_{g,n}(V)$ coincides with the moduli space of stable
curves $\M_{g,n}$.

There are two important canonical morphisms on $\M_{g,n}(V)$. The
{\em stabilization morphism\/}
\begin{equation} \label{eq:stabmor}
\st\,:\,\M_{g,n}(V)\,\to\,\M_{g,n}
\end{equation}
forgets the map $f$ and collapses the unstable components of the
curve $X$. The {\em evaluation morphism\/}
\begin{equation} \label{eq:evalmor}
\ev_i\,:\,\M_{g,n}(V)\,\to\,V
\end{equation}
 is obtained by evaluating the
stable map $f$ at the $i$-th marked point $p_i$.

The space $\M_{g,n}(V)$ has a special Chow homology class
\begin{equation}\label{eq:fundcl}
[\M_{g,n}(V)]^\mathrm{virt} \in H_*(\M_{g,n}(V),\mathbb{Q}),
\end{equation}
called the {\em virtual fundamental class\/}, that satisfies
various properties (axioms) formulated in~\cite{BMa}. The class
$[\M_{g,n}(V)]^\mathrm{virt}$ gives rise to the collection
$\Lambda^V\,\:=\,\{\,\Lambda^V_{g,n}\,\}$ of {\em Gromov-Witten
classes}
\begin{equation}\label{eq:gw}
 \Lambda^V_{g,n}\in H^\bullet(\M_{g,n})\,\otimes\,T^n H^\bullet(V)^*.
\end{equation}

We define $$
\Lambda^{V,\beta}_{g,n}(\gamma_1,\ldots,\gamma_n)\,:=\,\st_*
\left(
\ev_1^*\gamma_1\,\ldots\,\ev_n^*\gamma_n\,\cap\,[\M_{g,n}(V,\beta)]^\mathrm{virt}
\right), $$ where $\gamma_j\in H^\bullet(V)$.

\section{Stable Spin Maps and Their Moduli}

For the remainder of the paper we fix an integer $r\ge 2$. In this 
section we introduce the moduli space of stable $r$-spin maps and 
prove that it is a Deligne-Mumford stack.

\subsection{Overview of $r$-spin structures}

\

Although the concept of an $r$-spin structure is intuitively
simple, its formal definition is somewhat technical. For that
reason we first give a brief overview of the ideas involved.

\subsubsection{Spin structures on smooth curves}

\begin{df}\label{df:naive}
Given a collection of integers $\bm=(m_1, \dots, m_n)$, an
\emph{$r$-spin structure of type $\bm$} on a smooth, $n$-pointed
curve $(X, p_1, \dots, p_n)$ is a line bundle $\cl$ on $X$ with an
isomorphism $b:\cl^{\tensor r} \to \omega_X (-\sum m_ip_i)$.
\end{df}
In particular,  an $r$-spin structure of type
$\mathbf{0}=(0,\ldots,0)$ on a smooth curve $X$ is essentially
just an $r$-th root of the canonical bundle $\omega_X$.

For degree reasons,  an $r$-spin structure of type $\bm$ exists on
a genus $g$ curve $X$ only if $2g-2-\sum m_i$ is divisible by $r$.
When this condition is met, there are $r^{2g}$ choices of $\cl$ on
$X$.

\begin{df}
A \emph{smooth $r$-spin curve of type $\bm$} is a smooth curve $X$
with an $r$-spin structure of type $\bm$.  A \emph{smooth $r$-spin
map of type $\bm$ into a target $V$} is a morphism of a smooth
curve $X$ into $V$ with the additional data of an $r$-spin
structure  of type $\bm$ on $X$.
\end{df}

\begin{ex}
A $2$-spin structure of type $\mathbf{0}$ on a smooth curve $X$
corresponds to a choice of a theta characteristic $\cl$ on $X$ and
an explicit isomorphism $\cl^{\tensor 2} \to \omega_X$.
\end{ex}

\begin{ex}
If $E, p_0$ is a smooth elliptic curve, then $\omega_E$ is
isomorphic to $\co_E$.  An $r$-spin structure on $E$ of  type
$\mathbf{0}$ corresponds to a choice of an $r$-torsion point $q
\in E$ and an explicit isomorphism $$b: \co (q-p_0)^{\tensor r}
\irightarrow \co \cong \omega_E$$ from the $r$-th tensor power of
the invertible sheaf $\co(q-p_0)$ to the canonical sheaf
$\omega_E$.
\end{ex}

\subsubsection{Spin structures on nodal curves}

If we want to compactify the spaces involved by considering stable
maps and prestable curves, the above definition of an $r$-spin
structure is insufficient. In particular, even when the degree
condition is satisfied, there may be no line
 bundle $\cl$ on a prestable curve $X$ such that $\cl^{\tensor r}$ is
isomorphic to $\omega_X (-\sum m_ip_i)$.

The solution involves replacing line bundles by rank-one,
torsion-free sheaves, allowing the isomorphism $b:\cl^{\tensor r}
\to \omega_X (-\sum m_i p_i)$ to have non-trivial cokernel at the
nodes of the curve, and requiring that $b$ satisfy some additional
technical restrictions (see Definitions~\ref{df:root}
and~\ref{df:net}). There are two very different types of behavior
of this torsion-free sheaf $\cl$ near a node $q \in X$. When it is
still locally free, the sheaf $\cl$ is said to be \emph{Ramond} at
the node $q$. If the sheaf $\cl$ is not locally free at $q$, it is
called \emph{Neveu-Schwarz}.

In the Ramond case, the homomorphism $b$ is still an isomorphism
near the node $q$, but in the Neveu-Schwarz case it cannot be an
isomorphism. The local structure of the sheaf $\cl$ near a
Neveu-Schwarz node can be described as follows.

Near the node $q$, the curve $X$ has two coordinates $x$ and $y$,
such that $x y=0$.The sheaf $\omega_X$ (or $\omega_X (-\sum 
m_ip_i)$) is locally generated by $\frac{d x}{x} = -\frac{d 
y}{y}$.  Near $q$ the sheaf $\cl$ is generated by two elements 
$\ell_+$ and $\ell_-$, supported on the $x$ and $y$ branches 
respectively (that is, $x\ell_-=y\ell_+=0$).  The two generators 
may be chosen so that the homomorphism $b:\cl^{\tensor r} \to 
\omega_X(-\sum m_i p_i)$ takes $\ell^{\tensor r}_+$ to 
$x^{m_++1}(\frac{d x}{x})=x^{m_+}d x$     and $\ell^{\tensor 
r}_{-}$ to $y^{m_-+1}(\frac{d y}{y})=y^{m_-}d y$, where 
$(m_++1)+(m_-+1)=r$. 
\begin{df} We call $m_+$ (respectively
$m_-$) the \emph{order of the spin structure    along the
$x$-branch (respectively $y$-branch)}.
\end{df}

One more difficulty arises when $r$ is not prime---in this case
the moduli of stable curves with $r$-spin structure, as described
above, is not smooth.  The remedy is to include all $d$-spin
structures for every $d$ dividing $r$, satisfying some natural
compatibility conditions. This is described in Definition
\ref{df:net}.

\subsection{Definition of $r$-spin structures, curves, and maps}

\

We briefly review here the formal definition of an $r$-spin
structure over a fixed prestable curve, given in~\cite{J}, and we
define stable $r$-spin maps.

\begin{df}\label{df:root}
Let $(X, p_1, \dots, p_n)$ be a prestable, $n$-pointed, algebraic
curve; let $\ck$ be a rank-one, torsion-free sheaf on $X$; and let
$\bm = (m_1,\ldots,m_n)$  be an $n$-tuple of integers. A
\emph{$d$-th root of $\ck$ of type $\mathbf{m}$} is a pair $(\ce,
b)$, where $\ce$ is a rank-one, torsion-free sheaf, and $b$ is an
$\co_X$-module homomorphism $$  b: \ce^{\tensor d}  \rTo \ck
\tensor \co_X(-\sum m_ip_i) $$ with the following properties:
\begin{itemize}
\item $d \cdot \deg \ce = \deg \ck-\sum m_i$.
\item $b$ is an isomorphism on the locus of $X$ where $\ce$ is
locally free.
\item For every point $p \in X$ where $\ce$ is not free, the
length of the cokernel of $b$ at $p$ is $d-1$.
\end{itemize}
\end{df}

The condition on the cokernel amounts essentially to the condition
that at each node the sum of orders of $(\ce,b)$ is equal to
$d-2$.

If $\mathbf{m}'$ is congruent to $\mathbf{m}$ modulo $d$, then to
any $d$-th root $(\ce,b)$ of type $\mathbf{m}$ we can associate a
$d$-th root $(\ce',b')$ of type $\mathbf{m}'$ simply by taking
$\ce'=\ce \tensor \co(\frac{1}{d} \sum (m_i-m_i')p_i)$.
Consequently, the moduli of curves with $d$-th roots of a bundle
$\ck$ of type $\mathbf{m}$ is canonically isomorphic to the moduli
of curves with $d$-th roots of type $\mathbf{m}'$. Therefore,
unless otherwise stated, we will always assume the type
$\mathbf{m}$ of a $d$-th root satisfies $0 \leq m_i<d$ for all
$i$.

Unfortunately, when $d$ is not prime, the moduli space of curves
with $d$-th roots of a fixed sheaf $\ck$ is not smooth. To fix
this problem we must consider not just roots of a bundle, but
rather coherent nets of roots~\cite{J}.  This additional structure
suffices to make the moduli  space of curves with a coherent net
of roots smooth.

\begin{df}\label{df:net}
Let $\ck$ be a rank-one, torsion-free sheaf on a prestable,
$n$-pointed curve $(X, p_1, \ldots, p_n)$. A \emph{coherent net of
$r$-th roots of $\ck$ of type $\bm =(m_1, \ldots, m_n)$} is a pair
$(\{\ce_d\}, \{c_{d, d'}\})$ of a set of sheaves and a set of
homomorphisms as follows:  the set of sheaves consists of a 
rank-one, torsion-free sheaf $\ce_d$ on $X$ for every positive 
divisor $d$ of $r$; and the set of homomorphisms consists of an 
$\co_X$-module homomorphism $$ c_{d,d'} : \ce^{\tensor d/d'}_{d} 
\rTo \ce_{d'} $$ for every pair of divisors $d',d$ of $r$,  such 
that $d'$ divides $d$. These sheaves and homomorphisms must 
satisfy the following conditions: 
\begin{itemize}
\item $\ce_1=\ck$ and $c_{d,d}=\mathbf{1}_d$, the identity map, for every positive $d$ dividing $r$.
\item For each divisor $d$ of $r$ and each divisor $d'$ of $d$,
the  homomorphism $c_{d,d'}$  makes $(\ce_d, c_{d,d'})$ into a
$d/d'$-th root of $\ce_{d'}$ of type $\mathbf{m}'$, where
$\mathbf{m}'=(m'_1, \ldots, m_n')$ is the reduction of $\bm$
modulo $d/d'$ (i.e. $0\le m_i' < d/d'$ and $m_i \equiv m_i' \pmod
{d/d'}$).

\item The homomorphisms $\{c_{d,d'}\}$ are compatible.
That is, the diagram $$
\begin{diagram}
(\ce^{\tensor d/d'}_{d})^{\tensor d'/d''} &
\rTo^{(c_{d,d'})^{\tensor d'/d''}} & \qquad\ce^{\tensor
d'/d''}_{d'}\\ &  \rdTo^{c_{d,d''}}  & \dTo c_{d',d''} \\ & &
\ce_{d''}\\
\end{diagram}
$$ commutes for every $d''|d'|d|r$.
\end{itemize}
\end{df}

If $r$ is prime, then a coherent net of $r$-th roots is simply an
$r$-th root of $\ck$.  Even when $d$ is not prime, if the root
$\ce_d$ is locally free, then for every divisor $d'$ of $d$, the
sheaf $\ce_{d'}$ is uniquely determined up to an automorphism of
$\ce_{d'}$. In particular, if $\bm'$ satisfies the conditions
$\bm' \equiv \bm \pmod {d'}$ and $0 \leq m'_i <d'$,  then the
sheaf $\ce_{d'}$ is isomorphic to $\ce_{d}^{\tensor d/d'}\tensor
         \co\left(\frac{1}{d'} \sum (m_i -m'_i)p_i\right)$.

\begin{df}\label{def:spinstruct}
Let $X, p_1,\dots,p_n$ be an $n$-pointed, prestable curve of genus
$g$.  Let $r > 1$ be an integer and let $\bm = (m_1, m_2, \dots,
m_n)$ be an $n$-tuple of integers such that $r$ divides $2g-2-\sum
m_i$.  An \emph{$r$-spin structure on $X$ of type $\bm$} is a
coherent net of $r$-th roots of $\omega_X$ of type $\bm$ on $X$.
\end{df}

\begin{df}
Let $r > 1$ be an integer, and let $n$ and $g$ be non-negative
integers.  Let $V$ be a Deligne-Mumford stack, and let $\beta$ be
a class in $H_2 (V, \nz)$.  Finally, let $\bm = (m_1, m_2, \dots,
m_n)$ be an $n$-tuple of integers such that $r$ divides $2g-2-\sum
m_i$.  A \emph{stable, $n$-pointed, $r$-spin map into $V$ of genus
$g$, type $\bm$, and class $\beta$} is a pair
$(f,(\{\ce_d\},\{c_{d,d'}\}))$    that consists of a stable
$n$-pointed genus $g$ map  $f:X \to V$ of class $\beta$, and an
$r$-spin structure $(\{\ce_d\},\{c_{d,d'}\})$ of type $\bm$ on
$X$.
\end{df}

\begin{ex}
If $V$ is a point, then any stable, $n$-pointed, $r$-spin map into
$V$ is just a stable $r$-spin curve.
\end{ex}

\begin{df} \sloppypar
An isomorphism of $r$-spin maps $ (X \rTo^f V, p_1, \ldots, p_n,$
$(\{ \ce_d\}, \{c_{d, d'}\}))$ and $(X'\rTo^{f'} V, p'_1, \ldots,
p'_n, (\{\ce'_d\}, \{c'_{d,  d'}\})) $
 of the same type $\bm$ consists of
an isomorphism $\tau$ of $n$-pointed, stable maps $$
\begin{diagram}
X        & \rTo^{\tau} & X'\\ \dTo^{f} &             & \dTo^{f'}\\
V         & \rEq & V
\end{diagram}$$
and a set of $\co_X$-module isomorphisms $\{\theta_d : \tau^*
\ce'_d \irightarrow \ce_d\},$ with $\theta_1$ being the canonical
isomorphism $\tau^* \omega_{X'}(-\sum_i m_i {p'}_i) \irightarrow
\omega_X(-\sum m_ip_i),$ and such that the homomorphisms
$\theta_d$ are compatible with all the maps $c_{d,d'}$ and
$\tau^*c'_{d,d'}$.
\end{df}

Every $r$-spin structure on a smooth curve $X$ is determined, up
to isomorphism, by a choice of a line bundle $\ce_r$, such that
$\ce^{\tensor r}_r \cong \omega_X (-\sum m_i p_i)$. Therefore, in
the smooth case, the formal Definition~\ref{def:spinstruct} is
equivalent to the intuitive one from the previous subsection
(Definition~\ref{df:naive}). In particular, if $f:X\to V$ has no
automorphisms, then the set of isomorphism classes of $r$-spin
structures (if non-empty) of type $\mathbf{m}$ on $f:X\to V$ is a
principal homogeneous space for the group of $r$-torsion points of
the Jacobian of $X$. Thus there are $r^{2g}$ such isomorphism
classes.

\subsection{Families of $r$-spin structures and stable
spin maps} \label{sec:families}

\

To define the stack of stable $r$-spin maps, we must carefully
define how $r$-spin structures vary in families. This turns out to
be very delicate, since nilpotent elements may arise. In this
paper, we use the definition of families of spin curves given
in~\cite{J}.   The main condition is that all of the homomorphisms
$c_{d,d'}$ should be power maps in the sense of~\cite[\S
2.3.1]{J}. The definition for families given there reduces to the
definition of an $r$-spin structure given above, when the base is
(the spectrum of) a field; thus they are really only conditions on 
the families, rather than on the fibers.

The precise definition of families of spin curves is not necessary
for understanding the \cft\ related to the moduli space of stable 
$r$-spin maps discussed in 
Sections~\mbox{\ref{sec:classes}--\ref{sec:ex}}. In this paper we 
only use the formal definition in the proof of 
Theorem~\ref{thm:stab}, which states that the intuitive definition 
of the stabilization morphism from the stack of stable $r$-spin 
maps to the stack of stable $r$-spin curves does, indeed, give a 
morphism of stacks.  The reader willing to accept this result may 
skip the remainder of this subsection and the proof of 
Theorem~\ref{thm:stab} in Subsection~\ref{sec:proof}. 
\medskip

The following definition is rather technical, but it is precisely
what is needed to guarantee that the moduli stacks of spin curves
are smooth Deligne-Mumford stacks.  Families of spin curves can
also be defined (see~\cite{AJ}) in terms of line bundles on the
``twisted curves" of Abramovich and Vistoli~\cite{AV}.

\begin{df}\label{df:families}
An \emph{$r$-spin structure of type $\bm$} on a family  $X/T$ of
$n$-pointed prestable curves is a coherent net of $r$-th roots of
$\omega_{X/T}$ of type $\bm$.
\end{df}

Recall that, by Definition~2.3.4 of~\cite{J}, a {\em coherent net
of $r$-th roots\/} of $\omega_{X/T}$  of type $\bm$ is a set of
rank-one, torsion-free $\co_X$-modules $\{\ce_d\}$  and a
collection of $\co_X$-module homomorphisms $\{c_{d,d'}:
\ce_d^{\tensor d/d'} \to \ce_{d'}\}$, defined for  $d' | d | r$,
such that for each geometric fiber $X_t$ of $X/T$, the sheaves
$\{\ce_d\}$  and homomorphisms $\{c_{d,d'}\}$ induce a coherent
net of $r$-th roots of $\omega_{X_t}$ of type $\bm$, and each
homomorphism $c_{d,d'}$   is an isomorphism on the smooth locus of
$X/T$.  Finally, these sheaves and homomorphisms must have a
special type of local structure, which we describe now.

For a node $q\in X_t$ in a fiber of $X/T$ over a geometric point
$t\in T$, we denote by  $m_{d,+}$ and $m_{d,-}$ the orders of the
$d$-th root map $$ c_{d,1}: \ce_d^{\tensor d} \to \omega(-\sum m_i
p_i)$$ on the branches of the normalization of $X_t$ at $q$. We
define $$u_d:= (m_{d,+}+1)/\ell_d) \quad \mathrm{and} \quad  v_d:=
(m_{d,-}+1)/\ell_d), $$ where $$ \ell_d := \gcd(m_{d,+} + 1,
m_{d,-}+1). $$ If $c_{d,1}$ is an isomorphism at $q$, we set $u_d
= v_d = 0$.

The first requirement on the local structure of a net of coherent
roots on a family $X/T$  is the existence of a special local
coordinate system      near any node $q$ where $c_{r,1}$ is not an
isomorphism (i.e., $\ce_r$ is Neveu-Schwarz       at $q$). This
local coordinate system consists of an \'etale neighborhood $T'$
of $t$ with an element $\tau \in \co_{T',t}$, and an \'etale
neighborhood $U$ of $q$ in $X \cross_T T'$ with sections  $x,y \in
\co_{U}$, such that for $s:=u_r + v_r$ we have
   \begin{itemize}
\item $x y=\tau^{s}$.
\item The ideal generated by $x$ and $y$ has the singular
locus of $X/T$ as its associated closed subscheme.
\item The    homomorphism
$\big(\co_{T',t} [x,y]/(x y-\tau^{s})\big) \rightarrow \co_{U,q}$
induces an isomorphism of the completions $$
\left(\hat{\co}_{T',t} [[x,y]]/(x y-\tau^{s})\right) \irightarrow
\hat{\co}_{U,q}.$$
\end{itemize}

Of course, such a local coordinate system would always exist for
any $X/T$       if we had $s=1$, but in general    its existence
requires that $X/T$ be (locally) a ramified  cover of degree $s$
of another family of prestable curves.

The second requirement on the local structure is that the sheaves
$\ce_d$ must have a special presentation in terms of this special
coordinate system. In particular, any rank-one, torsion-free sheaf
$\cf$    always has a presentation of the form $$\cf \cong
\langle\zeta_1,\zeta_2 | e \zeta_1 = x \zeta_2, y \zeta_1 = h
\zeta_2\rangle $$ for some $e$ and $h$ in $\co_{T',t}$, such that
$e h = \tau^s$; but for sheaves in the net  we require that if
$\ce_d$ is not locally free at the node $q$, then $\ce_d$ must
have  such a presentation with $e=\tau^{(r/d)(v_d\ell_d)}$ and
$h=\tau^{(r/d)(u_d \ell_d)}$.

In other words, $\ce_d$ is isomorphic near node $q$ to the sheaf
$$ E_d:= \langle\zeta_1,\zeta_2 | \tau^{(r/d)(v_d \ell_d)} \zeta_1
= x \zeta_2, y \zeta_1 = \tau^{(r/d)(u_d\ell_d)} \zeta_2\rangle.$$
If $\ce_d$ is locally free at $q$, then for uniformity of notation
we will use the unusual presentation $\ce_d \cong E_d :=
\langle\zeta_1,\zeta_2| \zeta_1 = \zeta_2\rangle.$

Finally, each homomorphism $$c_{dj,j}: \ce_{dj}^{\tensor
d}\rTo\ce_{j}$$ in the net must be a so-called \emph{power map},
in the sense of Definition 2.3.1 of~\cite{J}. This means that, if
we use the local presentations $$ E_{dj} =\langle\xi_1,\xi_2 |
\tau^{(r/(dj))(v_{dj}\ell_{dj})} \xi_1 = x \xi_2, y \xi_1 =
\tau^{(r/(dj))(u_{dj} \ell_{dj})} \xi_2\rangle, $$ and $$ E_j =
\langle \zeta_1,\zeta_2 | \tau^{(r/j)(v_{j}\ell_{j})} \zeta_1 = x
\zeta_2, y \zeta_1 = \tau^{(r/j)(u_j \ell_{j})} \zeta_2\rangle,$$
of the sheaves $\ce_{dj}$ and $\ce_{j}$, then the map
\begin{equation}\label{eq:power}
\sym^d (E_{dj}) \to E_j,
\end{equation}
induced by the homomorphism $c_{dj,j}$, acts on the generators
$\xi_1^{d-i}\xi_2^{i}$ of  $\sym^d (E_{dj})$ as
            \begin{equation}\label{eq:powermap}
\xi_1^{d-i}\xi_2^{i} \mapsto \begin{cases} x^{u''-i}\tau^{iv}
\zeta_1 & \text{ if $0 \le i \le u''$} \\y^{v''- d + i}\tau^{(d-i)
v} \zeta_2 & \text{ if $u'' < i \le d$}. \end{cases}
\end{equation}
Here we require that $u_j \equiv u_{dj} d \pmod s$ and      define
$v_j \equiv v_{dj} d \pmod s$, \  $u'':= (u_{dj} d - u_j)/s$, and
$v'':= (v_{dj} d - v_j)/s$.

If $\ce_{dj}$ is locally free at $q$, then the existence of a good
presentation is automatically satisfied, and we have no additional
power map requirement except that the map~(\ref{eq:power}) be an
isomorphism.

This completes the definition of $r$-spin structures on families
of prestable curves, and we can now define families of $r$-spin
maps.

\begin{df}
Let $V$ be an algebraic variety and $\beta \in  H_2(V,\nz)$.  A 
\emph{family of  $n$-pointed, stable $r$-spin maps} to $V$ of 
class $\beta$ and type $\bm$ over a base $T$ is a family of 
$n$-pointed, stable maps $f:X \to V$ over $T$ of class $\beta$ 
with an $r$-spin structure of type $\bm$ on $X/T$. 
\end{df}

\subsection{Stacks of spin maps}

\

Now we are ready to define the main objects of the paper---the
stack of stable $r$-spin maps.

\begin{df}
Let $V$ be an algebraic variety over $\nc$, and $\beta$ an element 
of $H_2(V,\nz)$. The {\em stack of stable  $r$-spin  maps to 
$V$\/} ($n$-pointed, of genus $g$, and class $\beta$) is the 
disjoint union $$ \mgnrvbbar := \coprod_{\substack{\mathbf{m} \\ 0 
\leq m_i <r}}\mgnrmvbbar $$ of stacks $\mgnrmvbbar$ of (families 
of) stable $n$-pointed $r$-spin maps to $V$ of genus $g$, type 
$\bm=(m_1,\ldots,m_n)$, and class $\beta$. 
\end{df}

We will see in Section~\ref{sec:forget} that $\mgnrmvbbar$ (and,
therefore, $\mgnrvbbar$)  is a Deligne-Mumford stack whenever
$\mgnvbbar$ is. As the following proposition shows, no information
is lost by restricting $\bm$ to the range $0\le m_i \le r-1$.

\begin{prop}\label{prop:higherm}
If    $\bm \equiv  \bm' \pmod r$, then $\mgnrmvbbar$ is
canonically isomorphic to $\M^{1/r,\bm'}_{g,n}(V,\beta)$.
\end{prop}
\begin{proof}
When $\bm \equiv \bm' \pmod r$, every $r$-spin structure of type
$\bm$ naturally gives an $r$-spin structure of type   $\bm'$
simply by $$ \ce_d \mapsto \ce_d \tensor \co \left(\sum
\frac{m_i-{m'}_i}{d} p_i\right). $$
\end{proof}

\subsection{Canonical morphisms of stacks of stable  spin maps}

\

The stack $\mgnrmbar (V, \beta)$ has a natural projection
\begin{equation}\label{eq:forget}
\tilde{p}: \mgnrmbar (V, \beta) \rTo \mgnbar (V, \beta)
\end{equation}
which forgets the spin structure.  The usual evaluation maps
$$ev_i: \mgnbar (V, \beta) \to V,$$ which send a point $[X \rTo^f 
V, p_1 \dots p_n] \in \mgnbar (V, \beta)$ to $f(p_i) \in V$, 
induce evaluation maps $$\tilde{ev}_i =ev_i \circ \tilde{p} : 
\mgnrmbar (V, \beta) \rTo V. $$ 

Less obvious is the fact that for any morphism $s: V \to V'$
taking $\beta$ to $\beta'$, we have a stabilization morphism
\begin{equation}\label{eq:stab}
\tilde{st}: \mgnrmvbbar \rTo \mgnrmbar(V',\beta')
\end{equation}
which takes $f$ to $f' := s \circ f$ and contracts components of
the curve that are unstable with respect to $f'$.

\begin{thm}\label{thm:stab}
For any morphism $V \to V'$, taking $\beta$ to $\beta'$, the
stabilization map~(\ref{eq:stab}) is a morphism of stacks.
\end{thm}

The proof of Theorem~\ref{thm:stab} will be given in
Subsection~\ref{sec:proof}.
\medskip

The various canonical maps introduced above are shown in the
following commutative diagram.
\begin{equation}\label{eq:big-diagram}
\begin{diagram}
  &&                         & \mgnrmbar (V, \beta) &&                        && \\
  && \ldTo(3,4)^{\tilde{st}} & \dTo^{q_1} & \rdTo(3,4)^{\tilde{p}} \rdTo(5,4)^{\tilde{ev}_i} & \\
  &&                         & \mgnrmbar\kern-.75em \cross_{\mgnbar}\kern-.75em \mgnbar (V, \beta) &&&& & \\
  && \ldTo(3,2)^{pr_1}       &                      & \rdTo(3,2)^{pr_2}&  &        &   &\\
\mgnrmbar &&                       & \dTo^{q_2}           && &
\mgnbar (V, \beta) & \rTo^{ev_i} & V\\
  & \rdTo(3,2)_p            &&                      && \ldTo(3,2)_{st}              && &\\
  &&                         & \mgnbar              &&                         && &\\
\end{diagram}
\end{equation}
We will use the notation of this diagram throughout the remainder
of the paper, and we will denote the composition $q_2 \circ q_1$
by $q$.

The universal curves $\cgn \to \mgnbar$ and $\cgnrm \to
\mgnrmvbbar$ will be denoted by $\pi$.

\begin{rem}
The stack $\mgnrmvbbar$ is not isomorphic to the fibered product
$\mgnrmbar \cross_{\mgnbar} \mgnvbbar$, although on the smooth
locus the map $$q_1: \mgnrm (V, \beta) \rTo \mgnrm \cross_{\mgn}
\mgn (V, \beta) $$ is an isomorphism. This can be seen as follows.

If $X/T$ is a smooth family of curves, then a stable map $f: X \to
V$ and an $r$-spin structure $(\{\ce_d\}, \{c_{d,d'}\})$ are
precisely the data necessary to construct an $r$-spin map, i.e.,
there is a canonical morphism $$j: \mgnrm \cross_{\mgn} \mgn
(V,\beta) \to \mgnrm (V, \beta)$$ which is clearly the inverse of
the morphism $q_1$.

But if $X$ is not stable, this morphism $j$ no longer exists.  For
example, let $X$ be a prestable curve that has two irreducible
components $C$ and $E$, where $C$ is a smooth curve of genus $g$,
and $E$ is a smooth, rational curve, without marked points, joined
to $C$ at a single node $\fq$. Let $f: X \to V$ be an embedding of
$X$ in $V$.

An $r$-spin structure $(\{\ce_d\},\{c_{d,d'}\})$ on $X$ is
equivalent to a pair of  $r$-spin structures
$(\{\ce'_d\},\{c'_{d,d'}\})$ on $C$ and
$(\{\ce''_d\},\{c''_{d,d'}\})$ on $E$ of orders $0$  and $r-2$,
respectively, at $\fq$. Thus the automorphism group of the
$r$-spin map $(f, (\{\ce_d\}, \{c_{d,d'}\}))$ is $\mu_r \cross
\mu_r$, corresponding to multiplication of $\ce'_r$ and $\ce''_r$
by $r$-th roots of unity.

But the stabilization map $\tilde{st}$ takes $(f, (\{\ce_d\},
\{c_{d,d'}\}))$ to the spin map    $(f|_C,\{\ce'_d\},
\{c'_{d,d'}\})$ on $C$, and the automorphism group of
$$\tilde{p}(f, (\{\ce_d\}, \{c_{d,d'}\})) \cross \tilde{st}(f,
(\{\ce_d\},\{c_{d,d'}\}))=(f|_C, (\{\ce'_d\}), \{c'_{d,d'}\})$$ is
simply $\mu_r$, since $C$ is irreducible and $\ce'_r$ is
invertible on $C$.  Thus the morphism $q_1$ is not an isomorphism.
\end{rem}

\begin{prop}\label{q1-flat}
The morphism $q_1$ is flat and proper.
\end{prop}
\begin{proof}
Flatness follows from the valuative criterion of flatness
\cite[11.8.1]{EGA4}, which states that it is enough to check
flatness of $q_1$ over each $R$-valued point $\spec R \to
\mgnrmbar \cross_{\mgnbar} \mgnvbbar$, where $R$ is a discrete
valuation ring.  Since the completion $\hat R$ of $R$ is
faithfully flat over $R$, it suffices to check this for each
complete discrete valuation ring.  But in this case, the results
of \cite{J} show that the universal deformation (relative to the
universal stable map $f:\cc \to V$) of a spin structure over the
central fiber of $\spec R$  corresponds to the ring homomorphism
$R\to R[t]/(t^r-s)$, where $s\in R$ is a uniformizing parameter
for $R$.  In particular,  $R[t]/(t^r-s)$ is a free $R$-module, and
thus is flat over $R$.  Since the universal deformation is
faithfully flat (\'etale) over $\mgnrmvbbar$, this shows that
$q_1$ is also flat.

Properness also follows by the valuative criterion in exactly the 
same manner as was proved in \cite{J2} for spin structures on 
stable curves. Nothing in that proof required the underlying 
curves to be stable---only prestable.

\end{proof}

\subsection{The algebraic nature of the stack of stable spin maps}
\label{sec:forget}

\

A useful notion in dealing with stacks is the idea of a
{Deligne-Mumford morphism}, or {morphism of Deligne-Mumford type}.
This is analogous to the concept of a representable morphism.

\begin{df}
A morphism of stacks $f : S \to T$ is called {\em
Deligne-Mumford\/} (or \emph{of Deligne-Mumford type}) if for
every $U$-valued point $U \to T$, for a representable $U$, the
fibered product $S \cross_T U$ is a Deligne-Mumford stack.
\end{df}

The most useful fact about these morphisms is that if $S \to T $
is a Deligne-Mumford morphism, and if $T$ is a Deligne-Mumford
stack, then $S$ is a Deligne-Mumford stack (see~\cite[Prop.\
3.1.3]{J4}).

\begin{thm}\label{thm:forget}
For all $V$ and $\beta$, the forgetful morphism~(see
equation\ref{eq:forget}) is a finite Deligne-Mumford morphism of
stacks. In particular, $\mgnrmvbbar$ is a Deligne-Mumford stack
whenever $\mgnvbbar$ is.
\end{thm}

\begin{proof}
Given a $T$-valued point $T \to \mgnbar (V, \beta)$ for a
representable $T$, we must show that the stack $$R (X/T):=
\mgnrmbar (V, \beta) \cross_T \mgnbar (V, \beta)$$ of coherent
nets of $r$-th roots of $\omega_X (-\sum m_i p_i)$ on the
associated family $X/T$ of prestable curves is a Deligne-Mumford
stack, finite over $T$.  In particular, we need to construct a
smooth cover of $R(X/T)$ and show that the diagonal
$$\Delta:R(X/T) \cross_T R(X/T) \rTo R(X/T)$$ is representable,
unramified, and proper.

These facts are all straightforward generalizations of their
counterparts over the stack $\mgnbar$ of stable curves as
described in~\cite{J}.  The only real difference is that we are
now working with a specific family of prestable curves over $T$,
as opposed to working with the universal family of stable curves
(over $\mgnbar$), but that changes nothing of substance in the
proof.

The proof of properness is also an easy generalization of the case
of stable $r$-spin curves, and the morphism is obviously
quasi-finite, hence finite.
\end{proof}

\subsection{Proof of Theorem~\ref{thm:stab}}\label{sec:proof}

\

We now turn to the proof that for any morphism $s: V\to V'$,
taking $\beta \in H_2(V,\mathbb{Z})$ to $\beta' \in H_2
(V',\mathbb{Z})$, the  stabilization map~(\ref{eq:stab}) is a
morphism of stacks.

It is straightforward to check that the stabilization of the
underlying curves preserves $r$-spin structures on each individual
fiber, but we must also show that the stabilization morphism on
the underlying curves preserves the $r$-spin structure in
families.

Theorem~\ref{thm:stab} obviously follows from the following lemma.

\begin{lm}\label{lm:stab}
Let $st: \tilde{X}/T \to X/T$ be a morphism taking a family of
prestable curves $\tilde{X}/T$ to a partial stabilization $X$ of
$\tilde X$, and let $(\{\tilde{\ce_d}\}, \{\tilde{c}_{d,d'}\})$ be
an $r$-spin structure of type $\bm=(m_1,\dots,m_n)$ on
$\tilde{X}$, with $0 \le m_i \le r-1$ for every $i$.  In this
case, the sheaf $R^1st_*\tilde{\ce}_d$ is zero for every $d|r$,
and the push-forward $(\{st_* \tilde{\ce}_d\}, \{st_*
\tilde{c}_{d,d'}\})$ is an $r$-spin structure of type $\bm$ on
$X$.
\end{lm}

\begin{proof}
As mentioned above, it is straightforward to check that the maps
$st_* \tilde{c}_{d,d'}$ and the sheaves $st_* \tilde{\ce}_d$ are
$T$-flat and produce an $r$-spin structure of type $\bm$ on each
fiber of $X/T$ (this will also follow from the computations
below).  Thus we only need to verify that $R^1st_*\tilde{\ce}=0$
(which implies that it commutes with base change), and that the 
maps and sheaves meet the local conditions outlined in 
Subsection~\ref{sec:families} for being a coherent net on the 
family of curves $X/T$, provided the original sheaves 
$\{\tilde{\ce}_d \}$ and maps $\{\tilde{c}_{d,d'}\}$ form a 
coherent net on the family $\tilde X/T$. 

Let us fix  a point $p$ of a geometric fiber $X_t$ of $X/T$. There
are three cases to consider.  First is the case when the point $p$
is not the image of a contracted component (i.e., $st^{-1}(p)$ is
a single point).  Second is the case when $p$ is a smooth point of
the fiber $X_t$, but $p$ is the image of a whole irreducible
component of the fiber $\tilde{X}_t$ of $\tilde{X}/T$; that is,
$st$ contracts a $-1$-curve to the point $p$.  Third is the case
that $p$ is a node of the fiber $X_t$ containing it, and it is the
image of a contracted component of $\tilde{X}_t$; that is, $p$ is
the image of a $-2$-curve $\tilde{E}$.

\textbf{Case 1:} The first case is easy, since when $st^{-1}(p)$
is a single point, then $st$ is an isomorphism in a neighborhood
of $p$ (or of $st^{-1}(p)$).  In particular, $st_*$ is an
isomorphism,  $R^1st_*\tilde{\ce}_d=0$, and $c_{d,d'} = 
st_*(\tilde{c}_{d,d'})$ is a $d/d'$-th power map near $p$. 

\begin{figure}
\includegraphics{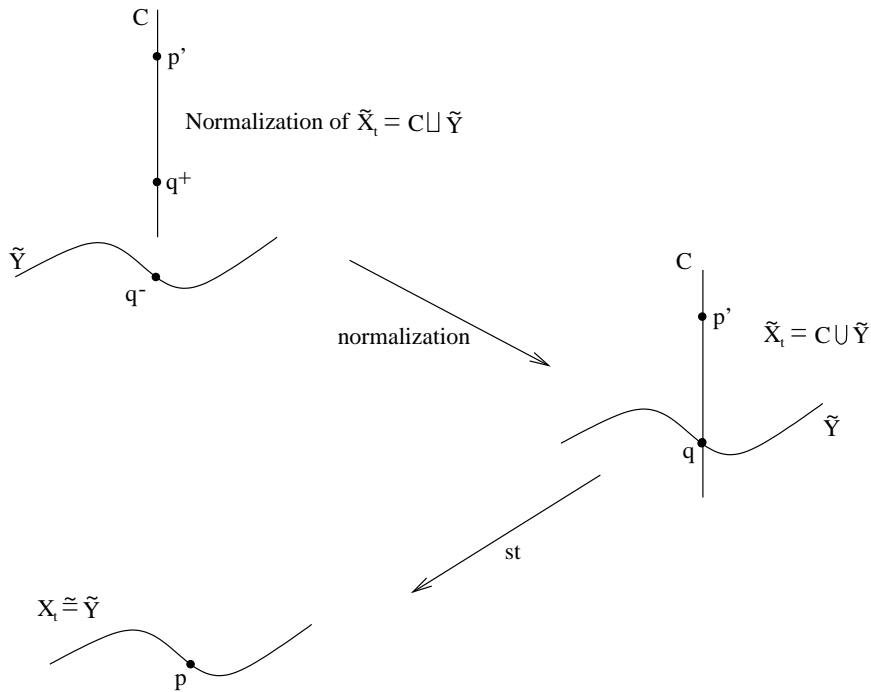}
\caption{\label{fg:case-2} A depiction of Case 2 of 
Lemma~\ref{lm:stab}:  fibers $\tilde X_t$ and $X_t$, the 
stabilization map $st: \tilde X_t \to X_t$, and the normalization 
of $\tilde X_t$. The morphism $st$ contracts the unstable 
component $C$ to the point $p$ and induces an isomorphism from the 
rest of the curve $\tilde Y$ to $X_t$.} 
\end{figure}

\medskip

The second and third cases are more involved.  Before we attack
them, we note that the conditions we must verify are local (and
analytic) on the base $T$, so it suffices to check the result when
$T$ is affine and is the spectrum of a complete local ring $R$.
Moreover, the conditions are analytic on $X$; that is, the
conditions are all determined by restricting to the completion of
the local ring of $X$ near the point $p$.  To simplify, we will
make the calculations in the case of $d=r$, but all other values
of $d$ (dividing $r$) are similar.

\textbf{Case 2:} In the second case ($st$ contracts a $-1$-curve
of $\tilde X$ to the point $p$) we will show that the induced
sheaves $st_* \tilde{\ce}_d$ are locally free at $p$, and the maps
$c_{d,d'}$ are all isomorphisms; thus the local coordinate and
power map conditions are automatically fulfilled.

The fiber $\tilde{X}_t$ over $X_t$ has one irreducible component
$C$ lying over $p$, and $C$ contains at most one marked point
$p'$, labeled with an integer $m$, where $0 \leq m \leq r-1$. This 
is indicated in Figure~\ref{fg:case-2}. On $C$, the sheaf 
$(\tilde{\ce}_r/\textsl{torsion})^{\tensor r}$ is isomorphic to 
$\omega_C(-m^+q^+-mp)$, where $q^+$ is the point of $C$ which maps 
to the node $q$ attaching $C$ to the rest of $\tilde{X}_t$.

Moreover, $r$ must divide $2g_C-2-m^+-m$, so either $m=r-1$, which
implies that $\tilde \ce_r$ is locally free (Ramond) near $q$, or          
$r-2=m^++m$, which implies that $\tilde \ce_r$ is not locally free 
(it is Neveu-Schwarz) at $q$.  In either case, $\tilde{\ce}_r|_C$ 
has degree $-1$ and thus has no global sections.  This also gives 
$R^1st_*\tilde \ce_r =0$, since this is true on each fiber.

Now, in the Neveu-Schwarz case, the sheaf
$st_*\tilde{\ce}_r|_{X_t}$ is simply the sheaf $\tilde{\ce}_r$
restricted (modulo torsion) to the rest of the prestable fiber
$\tilde{Y}=\overline{(\tilde{X}_t-C)}$.  But $
\tilde{\ce}_r/\textsl{torsion}$ on $\tilde{Y}$ is an $r$-th root
of $\omega_{\tilde{Y}}(-m^- q^- - \sum_{p_i \neq p} m_ip_i)$,
where $q^-$ is the other side of the node defined by $q^+$.  The
actual value of $m^-$ is determined by the relation $m^++m^-=m^+ +
m^- = r-2$, which implies that $m^-=m$.

In the Ramond case, the vanishing of the global sections of
$\tilde \ce_r|_C$ implies that $st_* (\tilde \ce_r|_{\tilde X_t})$
is $\tilde \ce_r|{\tilde Y} \tensor \ce(-q^-)$, so it is an $r$-th
root of $\omega_{X_t}(-(r-1)p)$.

In both the Ramond and Neveu-Schwarz cases, the new marked point
$p =st(q^-)$ of $X_t$ is labeled with $m$, just as the old marked
point $p'$ was marked with $m$ on $\tilde{X}_t$.  If no point was
marked on $C$, then the point $p$ remains unmarked (and $m^- =
0$).

Finally, $st_*\tilde{\ce}_d$ is $T$-flat and $R^1 st_* \tilde
\ce_r$ vanishes, so we have that $st_* \tilde \ce_r$ commutes with
base change, and the calculations above on the fibers all hold
globally on the family $X/T$.  Thus $st_*\ce_r$ is invertible near
$p$, and $st_*c_{r,1}$ is an isomorphism near $p$.  In particular,
$st_*c_{r,1}$ is an $r$-power map. A similar argument holds for
each $\tilde{\ce}_d$ and each $\tilde{c}_{d,d'}$ near $p$.

\textbf{Case 3:} The third case is that of a point $p \in X$ which
is the image of a $-2$-curve $\tilde{C}$ of $\tilde X$.  It is
easy to see that, just as in case 2,  on the unstable (contracted)
$-2$-curve, the degree of the bundle is $-1$.  Also, we have 
$R^1st_* \tilde \ce_r = 0$; the sheaf $st_* \tilde \ce_r$ is 
$T$-flat and commutes with base change; and on the fibers, the 
induced collection of sheaves and bundles forms an $r$-spin 
structure of type $\bm$.

We have still to check that the induced sheaves have the necessary
family structure for spin curves (existence of a local coordinate
of suitable type, with respect to which the sheaves have the
standard presentation---see Definition~\ref{df:families}), and
that the induced maps are power maps, as described in
equation~(\ref{eq:power}). For simplicity we will assume that the
orders $m_+$, $m_-$, $m'_+$, and $m'_-$ of the $r$-spin map
$\tilde{c}_{r,1}$ along the two nodes $q$ and $q'$ where the
$-2$-curve intersects the rest of the fiber have the property that
$\gcd(m_++1,m_-+1)=1 =\gcd(m'_++1,m'_-+1)$.  The case with common
divisors larger than $1$ is similar.

It is shown in~\cite[\S 3.1]{J2} that $\tilde{X}$ is locally
isomorphic to $$\proj_A A[\mu,\nu]/(x\nu-e^r\mu, h^r\nu-\mu y),$$
where $A =\hat{\co}_{X,x} \cong R[[x,y]]/xy-\pi^r$, and $e,h$ and
$\pi$ are elements of the maximal ideal $\maxid_R$ of $R$ with $e
h=\pi$.  This shows the existence of the special local coordinate.

We next show that $st_* \ce_{d}$ has a presentation of the form $$
st_* \ce_{d} \cong \langle\zeta_1,\zeta_2| \pi^{(r/d)(v_d \ell_d)} 
\zeta_1 = x \zeta_2, y  \zeta_1 = \pi^{(r/d)(u_d\ell_d)} 
\zeta_2\rangle.$$  If we let $\mu/\nu=s$ and $\nu/\mu=z$, then
near the exceptional $-2$-curve $\tilde{C}$ the curve $\tilde{X}$
is covered by two open sets, $$U = \{\mu \neq 0\} \cong \spec
A[z]/(x z-e^r,y-h^r z)$$ and $$ V=\{\nu \neq 0\} \cong \spec
A[s]/(x-e^rs, y s-h^r). $$

Since $(\{\tilde{\ce}_d\},\{c_{d,d'}\})$ is an $r$-spin structure,
we can describe $\tilde{\ce}_r$ on $U$ by $\tilde{\ce}_r|_U \cong
E_U(e^v,e^u) = \langle\zeta_1, \zeta_2|z \zeta_2 =e^{u}\zeta_1,
x\zeta_1=e^{v}\zeta_2\rangle$, and on $V$ by $\tilde{\ce}_r|_V
\cong E_V(h^{u'}, h^{v'})=\langle\xi_1, \xi_2|s \xi_2
=h^{u'}\xi_1, y \xi_1=h^{v'}\xi_2\rangle$, where $u+v=u'+v'=r$.

On the exceptional curve $\tilde{C}\cong \mathbb{P}^1$ the sheaf
$(\tilde{\ce}_r/\textsl{torsion})^{\tensor r}$ is isomorphic to
$\omega_{\mathbb{P}^1}((1-u)+(1-v'))$, and degree considerations
show that $u+v'=r$, so $u=u'$ and $v=v'$. Moreover, in a
neighborhood of $\tilde{C}$, if $D_i$ is the image of the $i$-th
section $\fp_i: T \to X$, the invertible sheaf $\omega_{\tilde{X}} 
(-\sum m_iD_i)$ is trivial and is generated by the element 
$w=\frac{dx}{x} =-\frac{dz}{z}=\frac{ds}{s}=-\frac{dy}{y}$.  The 
$r$-th power map $\tilde{c}_{r,1}$ is an isomorphism away from the 
nodes of $\tilde{X}$, and since it is a power map (changing the 
isomorphisms $\tilde{\ce}_r|_U\cong E_U(e^v,e^u)$ and 
$\tilde{\ce}_r|_V \cong E_V(h^u,h^v)$, if necessary), it maps the 
generators $\zeta_i$ and $\xi_i$ as follows: $$\zeta^r_i \mapsto 
z^u w, \,\, \zeta^r_2 \mapsto x^v w$$ and $$\xi^r_1 \mapsto s^vw, 
\,\,  \xi^r_2 \mapsto y^uw.$$

Since $\tilde{c}_{r,1}$ is an isomorphism away from the nodes, we
have $\zeta^r_1=z^r \xi^r_1$, or $\zeta_1=z\theta \xi_1$, for some
$r$-th root of unity $\theta$.  Changing the isomorphism
$\tilde{\ce}_r|_V \cong E_V(h^u,h^v)$ by $ \theta$, we may assume
$$\zeta_1=z\xi_1.$$

On $U \cap V$ we also have $$\zeta_2 =s e^v \zeta_1 =e^v\xi_1
\text{ and } \xi_2 =z h^u \xi_1=h^u\zeta_1.$$   So global sections 
of $\tilde{\ce}$ are of the form 
\begin{eqnarray*}
\Gamma(\tilde{\ce}_r)&=&\{((f_U\zeta_1+f'_U \zeta_2),
(f_V\xi_1+f'_V \xi_2))\in E_U \oplus E_V | \\
                 & &
f_U\zeta_1+f'_U\zeta_2=f_V\xi_1+f'_V\xi_2 \text{ on } U\cap V\}.
\end{eqnarray*}
We claim that the $A$-module
$$E(\pi^u,\pi^v):=\langle\eta_1,\eta_2|x\eta_2=\pi^u\eta_1,y\eta_1=\pi^v\eta_2\rangle
$$ is isomorphic to $\Gamma(\tilde{\ce}_r)$ via $$\eta_1 \mapsto
(\zeta_2, e^v \xi_1) \text{ and } \eta_2 \mapsto (h^u \zeta_1,
\xi_2).$$

The map is clearly an $A$-module homomorphism. Moreover, for any
section $((f_U \zeta_1+f'_U \zeta_2),(f_V \xi_1+f'_V\xi_2)) \in
\Gamma(\tilde{\ce}_r)$ we may assume that $f_U \in R[z]$ and $f'_U
\in R[[x]]$.  Likewise, we may assume that $f_V \in  R[s]$ and
$f'_V \in R[[y]].$

Consequently, we have $$z f_U(z) +e^v f'_U(x)-f_V(s)-z q^u
f'_V(y)=0,$$ or  $$z f_U(z)+e^v f'_U (s e^r)-f_V(s) - z h^u f'_V(z 
h^r)=0.$$ Thus $f_U$ and $f_V$ are completely determined by $$f_U 
=h^u f'_V(y) \text{ and } f_V=e^v f'_U(x).$$ 
 We may, therefore, map
$\Gamma(\tilde{\ce}_r)$ to $E(\pi^u, \pi^v)$ via $$(h^u
f'_V(y)\zeta_1+ f'_U(x)\zeta_2),(e^v f'_U(x)\xi_1+f'_V(y){\xi}_2)
\mapsto f'_U(x)\eta_1+f'_V(y)\eta_2,$$ and it is easy to check
that this homomorphism is the inverse of the first.

An identical argument shows that $\Gamma(\tilde{\ce}_d)$ is
isomorphic to $E(\pi^{u'}, \pi^{v'})$, where $u'\equiv u \pmod d$
and $v'\equiv v \pmod d$.  This shows the existence of the desired
presentation for $st_* \tilde{\ce}_d$.

It remains to show that the maps $st_* \tilde{c}_{d,d'}$ are power
maps~(\ref{eq:powermap}).  Again, since the arguments are
essentially identical for each pair $d$ and $d'$, it suffices to
prove this in the case of $\tilde{c}_{r,\sigma}$ for some $\sigma$
dividing $r$.

As above, we have $u+v=r$.  Let $\sigma$ be a divisor of $r$, and
$d=r/\sigma$. Let $u'$ be the smallest non-negative integer
congruent to $u d$ modulo $r$ and $v'$ be the smallest
non-negative integer congruent to $v d$ modulo $r$. Define
integers $u''$ and $v''$ as $$u''=\frac{d u-u'}{r} \ \text{ and }
\ v''=\frac{d v-v'}{r}.$$

The module $\Gamma(\tilde{\ce}_r) \cong E(\pi^u,\pi^v)$ is 
generated by $\eta_1,$ and $\eta_2$ with $\eta_1 =(\zeta_2, 
e^v\xi_1)$ and $\eta_2 =(h^u \zeta_1, \xi_2)$.  Further, 
$\tilde{\ce}_{\sigma}$ may be defined on $U$ by $\langle\phi_1, 
\phi_2 |z \phi_2 =e^{u'}\phi_1$, $x \phi_1 =e^{v'} \phi_2\rangle$ 
and on $V$ by $\langle\psi_1, \psi_2|s \psi_2=h^{v'}\psi_1, y 
\psi_1 =e^{u'}\psi_2\rangle$, so we may describe 
$st_*\tilde{\ce}_{\sigma}$ as above:  the module 
$\Gamma(\tilde{\ce}_{\sigma})$ is isomorphic to 
$E(\pi^{u'},\pi^{v'}),$ and is generated by $\gamma_1 =(\phi_2, 
e^{v'} \psi_1) $ and $\gamma_2 =(h^{u'} \phi_1, \psi_2)$.

We must show that $\eta^{d-i}_1 \eta^i_2$ maps, via
$st_*(\tilde{c}_{r,\sigma})$, to $\pi^{ui} x^{v''-i} \gamma_1$ for
$0 \leq i\leq u''$ and to $\pi^{v(d-i)} y^{u''-(d-i)}\gamma_2$
when $u''\leq i\leq d$.

We will do the first case---the second case is similar.  The
element $\eta^{d-i}_1 \eta^i_2$ is of the form $$\eta^{d-i}_1
\eta^i_2=(\zeta_2, e^v \xi_1)^{d-i} (h^u \zeta_1, \xi_2)^i=(h^{ui}
\zeta^i_1 \zeta^{d-i}_2, e^{v(d-i)}\xi^{d-i}_1 \xi^i_2),$$
 so on
$U$, this element $\eta^{d-i}_1 \eta^i_2$ maps as $$h^{ui}
\zeta^i_1 \zeta^{d-i}_2 \mapsto x^{v''-i}e^{ui}h^{ui} \phi_2
=\pi^{ui} x^{v''-i}\phi_2.$$
 On $V$, the element $\eta^{d-i}_1
\eta^i_2$ maps as $$e^{(d-i)v}\xi^{d-i}_1 \xi^i_2 \mapsto
s^{v''-i} h^{iu}e^{(d-i)v}\psi_1.$$

It is straightforward to check that these are the same on $U \cap
V$. But this is exactly the canonical $d$-th power
map~(\ref{eq:powermap}) for $E(\pi^v, \pi^u)^{\tensor d} \to
E(\pi^{v'}, \pi^{u'})$, as desired.
\end{proof}

\begin{rem}\label{rem:restrict-m} It is important to note that if any of the
$m_i$ is greater than $r-1$, Lemma~\ref{lm:stab} is no longer
true.  In particular, the sheaf $R^1st_* \tilde \ce_r$ no longer
vanishes in case~2 of the proof, and the subsequent
fiber-to-family transitions are not valid.
\end{rem}

\section{Cohomology Classes} \label{sec:classes}
\subsection{Tautological cohomology classes}

\

There are many natural cohomology classes in
$H^\bullet(\mgnrmvbbar,\nq)$; these include the classes induced by
pullback from $\mgnvbbar$ and from $\mgnrmbar$. In particular, we
have the $i$-th Chern  class $\lambda_i$ of the Hodge bundle
$\pi_* \omega_\pi$, and the components $\nu_i$ of the Chern
character of the Hodge bundle
\begin{equation}
\label{eq:nu} ch_t \pi_* \omega_\pi
    = 1+ch_t R\pi_* \omega_\pi
=g+\nu_1 t + \nu_2 t^3 +  \nu_3 t^5 + \dots~.
\end{equation}
(The even components of $ch_t \pi_* \omega_\pi$ vanish by
Mumford's theorem~\cite{Mu}.)  In a similar manner we define the
components $\mu_i$ of the Chern character of the pushforward
$R\pi_* \ce_r$ of the $r$-th root bundle $\ce_r$ to be
\begin{equation}
\label{eq:mu} ch_t R \pi_* \ce_r = -D+ \mu_1t + \mu_2 t^2 +
\dots~.
\end{equation}
Here, by $R \pi_* \ce_r$ we mean the K-theoretic pushforward
$\pi_* \ce_r - R^1 \pi_* \ce_r$, which is generally not the
equivalence class of a vector bundle, but only of a coherent
sheaf.\footnote{In \cite{JKV} the functor $R\pi_*$ was denoted by
$\pi_!$.} Here $-D$ is the Euler characteristic
$\chi(\ce_r|_{\cc_{{s}}})$ of $\ce_r$ on any geometric fiber
$\cc_{{s}}$ of $\pi$,  and by Riemann-Roch we have
\begin{equation}
\label{eq:degree} D=\frac{1}{r}\left((r-2)(g-1)+\sum_i m_i\right).
\end{equation}

In addition to the Hodge-like classes $\lambda_i$, $\nu_i$ and
$\mu_i$, there are tautological classes induced by the canonical
sections $$\fp_i:\mgnrmvbbar \to \cgnrm$$ corresponding to the 
marked points of $\pi: \cgnrm \to \mgnrmvbbar$. These are classes 
\begin{equation} \label{eq:psi}
\psi_i := c_1(\fp^*_i (\omega_{\pi})) \quad \mathrm{and} \quad
\tilde{\psi}_i := c_1(\fp^*_i (\ce_r))
\end{equation}
(and also class $\tilde{\psi}^{(d)}_i$ for each divisor $d$ of
$r$). When working in $\pic \mgnrmvbbar$, we will abuse notation
and use $\psi_i$ to indicate the line bundle
$\fp^*_i(\omega_{\pi})$, and $\tilde{\psi}_i$ the line bundle
$\fp^*_i(\ce_r)$. In~\cite{JKV} it is proved that these classes
are closely related:
\begin{equation} \tilde{\psi}_i= \frac{m_i+1}{r} \psi_i.
\end{equation}

\subsection{Graphs and boundary classes.}

\

Finally, there are the boundary classes. Much of the information
about the combinatorial structure of the boundary of $\mgnrmvbbar$
can be encoded in terms of decorated graphs.
\medskip

Recall that the \emph{(dual) graph} of an $n$-pointed, prestable
curve $(X, p_1,\ldots,p_n)$ consists of the following elements:
\begin{itemize}
\item \emph{Vertices},
corresponding to the irreducible components of $X$: a vertex $v$
is labeled with a non-negative integer $g(v)$, the (geometric)
genus of the component;
\item \emph{Edges},
 corresponding to the nodes of the curve:  an edge connects two
vertices (possibly even the same vertex, in which case the edge is
called a {\em loop}) if and only if the corresponding node lies on
the associated irreducible components;
\item \emph{Tails}, corresponding to the marked points $p_i \in X$,
$i=1,\ldots, n$: a tail labeled by the integer $i$ is attached at
the vertex associated to the component of $X$ that contains $p_i$.
\end{itemize}

\begin{df}
A \emph{half-edge} of a graph $\Gamma$ is either a tail or one of
the two ends of a ``real'' edge of $\Gamma$. We denote by
$V(\Gamma)$ the set of vertices of  $\Gamma$ and by $n(v)$ the
number of half-edges of $\Gamma$ at the vertex $v$.
\end{df}

\begin{df}
Let $\Gamma$ be a graph. The number $$ g(\Gamma) =
\mathrm{dim}H^1(\Gamma) + \sum_{v\in V(\Gamma)}g(v) $$ is called
the \emph{genus} of a graph $\Gamma$.
\end{df}

\begin{df}
A pair $(g,n)$ of non-negative integers is called \emph{stable} if
$2g+n-2 > 0$.

A vertex $v$ of a graph is called \emph{stable} if the pair
$(g(v),n(v))$ is stable.

A graph $\Gamma$ is called \emph{stable} if each vertex $v$ of
$\Gamma$ is stable.
\end{df}

To describe strata of the moduli space of stable $r$-spin maps
into $V$, we decorate the graphs with additional data coming from 
the $r$-spin structure and from $V$. In particular, the type 
$\mathbf{m} = (m_1,\ldots ,m_n)$ gives a marking to each of the 
tails, and the homology class $\beta_v$ of the image in $V$ of the 
corresponding irreducible component of the curve gives a marking 
to each vertex.

\begin{df}\label{df:stable-graph}
      Fix an integer $r\geq 2$ and a variety $V$. 
A \emph{($V,r$)-stable graph} is a graph $\Gamma$ with a choice of 
a homology class $\beta_v \in H_2(V,\nz)$ for each vertex $v$ of 
$\Gamma$, and a marking of each half edge $h$ by a non-negative 
integer $m_h<r$. For each edge $e$ the marks $m_+:=m_{h_+}$ and 
$m_-:=m_{h_-}$ of the two half-edges $h_+$ and $h_-$ of $e$ must 
satisfy 
\begin{equation}  \label{eq:edgecong}
 m_+ + m_- \equiv r-2 \pmod r.
\end{equation}

Finally, the graph should satisfy the stability condition that if
$\beta_v=0,$ then the vertex $v$ is stable. 

If the  half-edges of $\Gamma$ are not marked by integers $m_h$,
but its vertices are marked with classes $\beta_v$, and all
vertices with $\beta_v=0$ are stable, then such a graph will be
called a \emph{$V$-stable graph}.

Similarly, if the vertices of $\Gamma$ are stable and not marked,
but the half-edges are marked with integers $m_h$
satisfying~(\ref{eq:edgecong}), then it will be called an
\emph{$r$-stable graph}.
\end{df}

Each stable $r$-spin map defines a $(V,r)$-stable graph, called
its dual graph.

\begin{df}
Given a stable $r$-spin map $f:X\rTo V, (\{\ce_d\},\{c_{d,d'}\})$,
its  \emph{$(V,r)$-decorated dual graph}  (or just \emph{dual
graph}) is the dual graph $\Gamma$ of the underlying curve $X$,
with the following additional markings.  Each vertex $v$ is
labeled with the class $\beta_v := [f(X_v)]$ of the image of the
corresponding irreducible component $X_v$.  The $i$-th tail is
marked by $m_i$, and each half-edge associated to a node of $X$ is
marked by the order (the integer $m^+$ or $m^-$) of the $r$-spin
structure along the branch of the node associated to that
half-edge.
\end{df}

For any $V$, we let $G^{1/r}_{g,n}(V)$ denote the set of all
$(V,r)$-stable  graphs of genus $g$ with $n$ tails.

\begin{df-pr}
For any morphism $V \rTo^{\gamma} V'$ and a stable pair $(g,n)$, 
there is an associated stabilization map 
\begin{equation} \label{eq:graphstab}
\gamma_{\bullet}:G^{1/r}_{g,n}(V)\to G^{1/r}_{g,n}(V'),
\end{equation}
defined as follows. The graph $\Gamma \in G^{1/r}_{g,n}(V)$ is
mapped to the graph $\overline{\Gamma} \in G^{1/r}_{g,n}(V')$ 
obtained from $\Gamma$ by removing all vertices $v$ in $\Gamma$ 
that fail the stability criterion; that is, $\gamma_*(\beta_V)=0$ 
in $H_2(V',\mathbb{Z})$ and  $$2g(v)-2 +n(v)\le 0.$$

Also remove all half-edges attached to each removed vertex $v$,
and join together any other half-edges that were previously
connected to half-edges of $v$. Since there are at most two such
half-edges per unstable vertex $v$, this operation will either
produce a well-defined edge, or it will not connect any half-edges
at all.  Now mark each remaining vertex $v$ of the resulting graph
with $\gamma_*(\beta_v)$, and give each half-edge $h$ the mark
$m_h$ it had in the graph $\Gamma$.  The resulting graph
$\overline{\Gamma}$ is clearly $(V',r)$-stable, except in the
special case that for every vertex $v$ of $\Gamma$ we have
$\gamma_* (\beta_v)=0$ and $2g(v)-2+n(v)\leq 0$.  This only occurs
when $2g-2+n\leq 0$, so when $2g-2+n>0$ we always have a
stabilization map $\gamma_{\bullet}:G^{1/r}_{g,n}(V) \to
G^{1/r}_{g,n}(V')$ associated to $\gamma$.
\end{df-pr}

If the morphism $\gamma$ is the constant map taking $V$ to a
point,
 we will call the image of a $V$-stable (or $(V,r)$-stable) graph
$\Gamma$ under the associated stabilization map
$\gamma_{\bullet}$ just \emph{the stabilization of $\Gamma$}.
\medskip

$(V,r)$-stable graphs of genus $g$ with $n$ tails correspond to
boundary strata in $\mgnrvbbar$, although some of these strata may
be empty.

\begin{df}
Let $\Gamma$ be a connected $V$-stable graph (or $(V,r)$-stable,
or $r$-stable graph) with $n$ tails and of genus $g$. We denote by
$\mgammavbbar$ (or by $\mgammarvbbar$, or by $\mgammarbar$) the
closure in $\mgnvbbar$ (or in $\mgnrvbbar$, or in $\mgnrbar$) of
the moduli space of stable maps (or stable $r$-spin maps, or
stable $r$-spin curves) whose dual graph is $\Gamma$. If
$\Gamma=\amalg \Gamma_i$ is the disjoint union of connected
subgraphs $\Gamma_i$, then we denote by $\M_{\Gamma}(V)$ the
product $\amalg \,\M_{\Gamma_i}(V)$, and similarly
$\M_{\Gamma}^{1/r}(V) = \amalg \,\M^{1/r}_{\Gamma_i}(V)$.
\end{df}

These substacks $\mgammarvbbar$ are generally not irreducible,
since different spin structures often cannot be deformed into one
another, even for a fixed type $\bm$. For example, when $r = 2$,
$g > 0$, $\bm = \mathbf 0$,  and the target $V$ is a point, there
are both even and odd spin structures on each stable curve, and
these form distinct irreducible components.

Moreover, the classes $[\mgammarvbbar]$ in the Chow group
$A_*(\mgnrmvbbar)$ defined by the substacks $\mgammarvbbar$ are
not usually the pullbacks of the corresponding classes in
$A_*(\mgnvbbar)$.  For example, if $V$ is a point and the graph
$\Gamma$ is a tree with one edge and two vertices $$\Gamma
=\tcongm,$$ then there is a unique choice of marking $m^+$ and
$m^-$ on the two half-edges of the edge that makes the degree of
the twisted canonical bundle divisible by $r$ on both vertices. In
this case it can be shown (see~\cite{J3}) that in $\pic \mgnrmbar$
the class $\tilde{p}^* [\mgammabar]$ is precisely $$ \tilde{p}^*
[\mgammabar] = \frac{r}{\gcd(m^++1,r)} [\mgammarbar].$$

\section{Spin Virtual Class}  \label{sec:virt}

Recall from~\cite[\S4.1]{JKV} that an $r$-spin virtual class on
the moduli of stable, $r$-spin curves is an assignment of a
cohomology class \begin{equation} \label{eq:cvirt} \cv_\Gamma \in
H^{2 D} (\mbar^{1/r}_\Gamma, \nq)
\end{equation}
 to every genus $g$, $r$-stable graph $\Gamma$ with
$n$-tails. Here, if the tails of $\Gamma$ are marked with the
$n$-tuple $\bm=(m_1,\ldots,m_n)$, then the dimension $D$ is
\begin{equation}
\label{eq:deg} D = \frac{1}{r}\bigl( (r-2)(g-\alpha)+\sum_{i=1}^n
m_i \bigr),
\end{equation}
where $\alpha$ is the number of connected components of $\Gamma$.
In the special case where $\Gamma$ has one vertex and no edges, we
denote $\cv_\Gamma$ by $\cv_{g,n}(\bm)$.  The classes are required
to satisfy the axioms of connected and disconnected graphs,
convexity, cutting edges, vanishing, and forgetting tails.

We can use the choice of an $r$-spin virtual class for stable
$r$-spin curves to produce a similar $r$-spin class for all stable
map spaces.

\begin{df}\label{class} Let $G^{1/r}$ be the set of all $r$-stable graphs.
Given an $r$-spin virtual class $\{\cv_{\Gamma} \in
H^{2D}(\M^{1/r}_{\Gamma}, \nq)\}_{\Gamma \in G^{1/r}}$ meeting the
axioms of~\cite[\S4.1]{JKV}, then for each $V$ and for each
$(V,r)$-stable graph $\Gamma$, define the class
\begin{equation}\label{eq:def1}
\cvt_{\Gamma}= \tilde{st}^* \cv_{\overline{\Gamma}},
\end{equation}
where $\overline{\Gamma}$ is the stabilization of $\Gamma$.
\end{df}

The only graphs that are $V$-stable but have no stabilization are
graphs with $2g-2+n\leq 0$, so in these cases we define the
$r$-spin virtual class directly.

\begin{df}\label{chern}
If $g=0$ and $n<3$ then we define $\cvt_{\Gamma}$ to be the top
Chern class of the dual of the first cohomology of the $r$-th root
bundle $\ce_r$; namely,
\begin{equation}\label{eq:def2}
\cvt_{\Gamma}=c_D (-R^1 \pi_* \ce_r),
\end{equation}
where $\ce_r$ is the $r$-th root of the universal spin structure
$(\{\ce_d\},\{c_{d,d'}\})$ on the universal curve $\pi: \cc \to
\M^{1/r}_{\Gamma} (V, \beta).$
\end{df}

\begin{prop}\label{pr:unstab}
If $g=0$ and $n<3$, then for every connected graph $\Gamma \in
G^{1/r}_{0,n}(V)$ with no marking $m_i$ equal to $r-1$, the
$r$-spin virtual class $\cvt_{\Gamma}$ has dimension zero; and
thus for any $\Gamma \in   G^{1/r}_{0,n}(V)$ with $n<3$ we have
\[
\cvt_{\Gamma}=\begin{cases} 0 & \text{if any } m_i =r-1
\\ 1 & \text{otherwise}.  \end{cases}    \]
\end{prop}

\begin{proof}
The degree of the sheaf $\ce_r$ is an integer and is given by
$$\deg \ce_r =(2g-2-\sum m_i)/r,$$ hence when $g=0$ we have $$\sum
m_i \equiv -2 \pmod r.$$  The dimension $D$ of $\cvt_{\Gamma}$ is
$$D = ((2-r) +\sum m_i)/r.$$  If $n=0$ we have $\sum m_i=0$, which
implies $r=2$, and we immediately have $D=0$.  If $2 \geq n \geq
1$, then since $0 \leq m_i \leq r-2$, we have $0 \leq \sum m_i 
\leq 2(r-2)$; and hence $\sum m_i=r-2$ is the only solution to the 
congruence $\sum m_i \equiv -2 \pmod r$. Consequently, 
$$D=(2-r+\sum m_i)/r=0.$$

If any of the $m_i$ are equal to $r-1$, then the argument in the
proof of Axiom 4 in~\cite[Theorem~4.1]{JKV} shows that $\cvt$ must
be zero.
\end{proof}

In the case that $g=1$ and $n=0$, the moduli space $\M^{1/r}_{1,0}
(V, \beta)$ decomposes into the disjoint union of $d$ connected
components, where $d$ is the number of positive divisors of $r$
(including $1$ and $r$); these components correspond to the fact
that (on the smooth locus) $r$-spin structures are in one-to-one
correspondence with $r$-torsion points of the Jacobian of the
underlying curve.  No deformation of the underlying curve can take
a point of order $i$ to a point of order $j$ unless $i=j$, so the
moduli space breaks up into disjoint components $$\M^{1/r}_{1,0}
(V, \beta) = \coprod_{\substack{i|r \\ 1 \leq i\leq r}}
\M^{1/r,(i)}_{1,0} (V, \beta).$$ We call $i$ the \emph{index} of
the component if the $r$-th root is a point of exact order $i$ in
the Jacobian of the underlying curve.

\begin{df}\label{r-spin}
If $g=1$ and $n=0$, define the $r$-spin virtual class
$\cvt_{\Gamma} (V, \beta)$ as follows
\begin{equation}\label{eq:def3}
\cvt_{\Gamma} (V, \beta) =\left\{ \begin{array}{ll} -(r-1) &
\text{if the index is $1$}\\ 1 & \text{otherwise.}
\end{array} \right.
\end{equation}
\end{df}

\begin{thm}\label{thm:cv-convex}
If $g=0$, then $\cvt_{\Gamma} (V, \beta)= \tilde{st}^*
\cv_{\overline{\Gamma}}$ is the top Chern class $c_D (-R^1 
\tilde{\pi}_* \tilde{\ce}_r)$ of the bundle whose fiber is the 
dual of the first cohomology of the $r$-th root $\tilde{\ce}_r$ on 
the universal curve $\tilde{\pi}: \tilde{\cc} \to 
\M^{1/r}_{\Gamma}(V, \beta)$. 
\end{thm}

\begin{proof}
For $n < 3$, this is true by definition.

In the case that $n\ge 3$, since $g=0$, the $r$-spin virtual class
$\cv_{\Gamma} \in H^{2D}(\M^{1/r}_{\Gamma}, \nq)$ is the  top
Chern class $c_D (-R^1 \pi_* \ce_r)$ of the first cohomology of
the $r$-th root $\ce_r$ on the universal curve $\pi: \cc\to
\M^{1/r}_{\overline{\Gamma}}$, by the convexity axiom of
\cite{JKV} \S 4.1.

We have the following commutative diagram. $$\begin{diagram}
\tilde{\cc} & \rTo^{\phi} & \cc
\cross_{\M^{1/r}_{\overline{\Gamma}}} \M^{1/r}_{\Gamma} (V, \beta)
& \rTo^{p_1} & \cc \\ & \rdTo_{\tilde{\pi}} & \dTo_{p_2} & &
\dTo_{\pi} \\ & & \M^{1/r}_{\Gamma} (V, \beta) & \rTo^{\tilde{st}}
& \M^{1/r}_{\overline{\Gamma}} \end{diagram}.$$

Here $\phi$ is the natural map induced by $\tilde{\pi}$ and
stabilization of $\tilde{\cc}$.  If $\tilde{\ce}_r$ is the $r$-th
root on $\tilde{\cc}$, then by Lemma \ref{lm:stab} and the
universality of the sheaves involved, $\tilde{\phi}_*
\tilde{\ce}_r$ is isomorphic to the pullback $p^*_1 \ce_r$ of the
$r$-th root $\ce_r$ from $\cc$, and $R^1 \phi_* \tilde{\ce}_r =0$.
By the Leray spectral sequence we have $$R^1 \tilde{\pi}_*
\tilde{\ce}_r = R^1 p_{2*} (p^*_1 \ce_r).$$ But $\tilde{st}$ is
flat (it is the composition of flat morphisms---see the
commutative diagram~(\ref{eq:big-diagram})), so that $$c_D(-R^1
p_{2*}(p^*_1(p^*_1 \ce_r)) =\tilde{st}^* c_D (-R^1 \pi_*
\ce_r)=\tilde{st}^* c^{1/r}_{\Gamma}.$$
\end{proof}

\begin{rem}\label{rem:cv-m-restrict}
The proof of Theorem~\ref{thm:cv-convex} depends upon the fact
that the integer marking $m_h$ of each half edge $h$ lie in the
range $0 \le m_h \le r-1$, as required for stable graphs (see
Definition~\ref{df:stable-graph}).   In particular, when an $m_h$
lies outside that range, Lemma~\ref{lm:stab} fails.

We shall see (in Remark~\ref{rem:descent}) that
Theorem~\ref{thm:cv-convex} is false in the case that any $m_h$ is
larger than $r-1$.
\end{rem}

\begin{df}
We define $[\mgnrvbbar]^{\mathrm virt}$ to be the pullback
$$[\mgnrvbbar]^{\mathrm virt} := \tilde p^* [\mgnvbbar]^{ \mathrm
virt}$$ of the usual virtual fundamental class
$[\mgnvbbar]^{\mathrm virt}$ (see Section~\ref{sec:fundcl}) of
$\mgnvbbar$ via $$\tilde p: \mgnrvbbar \to \mgnvbbar$$
\end{df}

Using the notation of the commutative
diagram~(\ref{eq:big-diagram}), since $\evt_i =\ev_i \circ
\tilde{p}$, for any $\gamma_1, \dots, \gamma_n \in H^{\bullet} (V,
\nq)$ we have the equality
     $$ \evt_1^* \gamma_1 \cup
\evt^*_2 \gamma_2 \cup \dots \cup \evt^*_n \gamma_n
=\tilde{p}^*(\ev^*_1 \gamma_1 \cup \dots \cup \ev^*_n \gamma_n). $$

We also have the following important relation on classes, which is
the main step in proving that the $\cft$ defined by stable
$r$-spin maps is the tensor product of the \cfts\ of $r$-spin
curves and stable maps (Theorem~\ref{thm:tensorproduct}).

\begin{thm}
\label{thm:tensor} Given any set $\{\gamma_1, \dots, \gamma_n\}$
of classes in $A^*(V)$ (or $H^\bullet(V)$), and given an $r$-spin 
virtual class $\cvt$ on $\mgnrmbar (V, \beta)$ defined by 
equations~(\ref{eq:def1}),~(\ref{eq:def2}) and ~(\ref{eq:def3}), 
the relation 
\begin{equation}\label{eq:product}
q_* (\cvt \cup \prod^{n}_{i=1} \evt^*_i (\gamma_i)\cap
[\M_{g,n}^{1/r}(V,\beta)]^\mathrm{virt})=p_* \cv \cup st_*
(\prod^n_{i=1} \ev^*_i (\gamma) \cap [\M_{g,n}(V,\beta)
]^\mathrm{virt})
\end{equation}
holds.
\end{thm}
\begin{proof}
We will give the proof on the level of (operational) Chow groups
$A^*$ with notation as in~\cite[V \S8]{Ma}.  From \cite[VI 
\S2]{Ma} it will follow then that such results also hold for 
$H^\bullet(V)$.

In particular, if we denote the identity maps on $\mgnbar$,
$\mgnrbar$, $\mgnvbbar$, $\mgnrbar \cross_{\mgnbar} \mgnvbbar$,
and $\mgnrvbbar$ by $\mathbb{I}$, $\mathbb{I}_r$,
$\mathbb{I}_{V}$, $\mathbb{I}_{\cross}$, and $\mathbb{I}_{r,V}$,
respectively, then we have $\cv \in A^*(\mgnrbar):=
\bar{A}^*(\mathbb{I}_r: \mgnrbar \to \mgnrbar)$, and $\cvt =
\tilde{st}^*(\cv)\in A^*(\mgnrvbbar)$. We take $\gamma_i$ in
$A^*(V) $, so that $\evt_i^*(\gamma_i)$ is in $ A^*(\mgnrvbbar)$.
Also, we have $[\mgnvbbar]^{\mathrm{virt}} \in A_*(\mgnvbbar)$.

Finally,  by $$\cvt \cup \prod^n_{i=1} \evt_i^*(\gamma_i) \cap
[\M_{g,n}^{1/r}(V,\beta)]^{\mathrm{virt}}$$ we mean $$\left(\cvt 
\cup \prod^n_{i=1} \evt_i(\gamma_i)\right)_{\mathbb{I}_{r,V}} \cap 
[\M_{g,n}^{1/r}(V,\beta)]^{\mathrm{virt}}.$$ 

As in~\cite[V \S8.9]{Ma}, for any morphism $Y \to X$, we define
$f^*:A^*(X) \to A^* (Y)$ to be \begin{equation} f^* (\delta)_h
\cap y:=\delta_{f \circ h} \cap y, \end{equation}where $\delta \in
A^* (X)$ and $h: L \to Y$ is an arbitrary morphism, and $y \in A_*
(L)$.  We also define, for any proper, flat morphism $f:Y\to X$ of
Deligne-Mumford stacks $X$ and $Y$, the proper flat pushforward
$f_{\bullet}: A^*(Y) \to A^*(X)$ to be
\begin{equation} f_{\bullet} \alpha_g \cap c := f_*(\alpha_{f'_Y} \cap
f^*(c)),\end{equation} where $g:L\to X$ is an arbitrary morphism,
$\alpha$ is an element of $A^*(Y)$, and $c$ is an element of
$A_*(X)$.

\begin{rem}
Note that part (ii) of Manin's definition in \cite[V \S 8.9]{Ma}
of the operational Chow ring $A^*(M)$ for the identity morphism
$\mathbb I: M \to M$ states that elements of  $A^*(M)$ only need
to commute with pullback along \emph{representable}, flat
morphisms of DM-stacks, despite the fact that standard definitions
of general operational Chow rings require that these elements
commute with pullback along \emph{all} flat morphisms of DM-stacks 
(see Vistoli \cite[5.1.i]{V} and Manin \cite[V.8.1.i]{Ma}).

In what we do below, we will need the definition of $A^*(M)$ that
requires commutativity with all flat pullbacks; that is, we
require the following.

Let $f:X \to Y$ be a flat
 morphism of Deligne-Mumford stacks, which is \emph{not necessarily
representable}, and let $h: Y \to Z$ be an arbitrary morphism of
Deligne-Mumford stacks.  For any $\sigma \in A^{*}(Z)$ and $y \in
A_*(Y)$, we have
\begin{equation}\label{eqn:flat-pb} \sigma_{h \circ f} \cap f^*
(y) =f^*(\sigma_h \cap y).\end{equation}

This seemingly minor difference in the definition of $A^*$ allows
us to prove a projection formula for non-representable morphisms.

\end{rem}

\begin{lm}\label{flat}

\

\begin{enumerate}
\item[1.] Let $f:X\to Y$ be a proper, flat morphism of
Deligne-Mumford stacks (which is not necessarily representable),
and let $h:L \to Y$ be an arbitrary morphism of Deligne-Mumford
stacks.  We have
\begin{equation}\label{eqn:flat-pb2} h^* f_{\bullet} =
f_{L\bullet} h^*_X.\end{equation}

\item[2.] \textbf{(Projection formula for $f_{\bullet}$)} Let $f:X \to Y$ be a proper, flat morphism of
Deligne-Mumford stacks (which is not necessarily representable).
For any $\sigma \in A^* (X)$ and $\beta \in A^*(Y)$ we have
\begin{equation}
f_{\bullet} (\sigma f^* (\beta)) =f_{\bullet} (\sigma) \beta.
\end{equation}
\end{enumerate}
\end{lm}

\begin{proof}
For  part 1 of the lemma, the same proof as given by Manin for
this equation  \cite[V.8.30]{Ma} works exactly for our case, too; 
nowhere is the representability of $f$ used in Manin's proof.

For part 2, again Manin's proof of the projection formula
\cite[V.8.29]{Ma} works for non-representable morphisms, the only
change needed is that~\cite[V.8.22]{Ma} (commutativity with flat,
representable pullbacks) must be replaced by our equation
(\ref{eqn:flat-pb}) for non-representable, flat pullbacks.
\end{proof}

One more fact we will need in the proof of
Theorem~\ref{thm:tensor} is the commutativity with proper
pushforwards required by the definition of $A^*$ (cf.\
\cite[V.8.21]{Ma}); namely, if $p:P \to L$ is proper, and $h:L\to 
M$ is an arbitrary morphism, then by definition of $A^*(M)$, for 
any $\sigma \in A^*(M)$ and for any $y \in A_*(P)$ we have 
\begin{equation}\label{eqn:proper-push}
\sigma_h \cap p_* (y) =p_*(\sigma_{hp} \cap y).
\end{equation}

Now we may proceed with the proof of Theorem~\ref{thm:tensor}.  We
will refer throughout the proof to the notation of the commutative
diagram (\ref{eq:big-diagram}).

Since $q_1$ is a birational map, it is a splitting morphism (i.e.,
$q_{1\bullet} q_1^* = \mathbb{I}_{\cross}$, as a map on $A^*$).
Moreover, the morphism $st$ is flat~\cite[Prop. 3]{B2} and proper,
$p$ is flat \cite[Theorem 2.2]{J} and proper \cite[Theorem
2.3]{J}, and $q_1$ is flat and proper by Proposition 
\ref{q1-flat}. We have the following relations: 

$$q_* \left( \left(\cvt \cup \prod^n_{i=1}
\tilde{ev}_i(\gamma_i)\right) \cap
\tilde{p}^*[\mgnvbbar]^{\mathrm{virt}}\right) $$
\begin{eqnarray*}
&  =&  q_{2*} q_{1*}\left ( q^*_1 \left(pr^*_1 \cv \cup pr^*_2
\prod^n_{i-1}ev^*_i (\gamma_i)\right)_{\mathbb{I}_{r,V}} \cap
q_{1}^* pr_2^*[\mgnvbbar]^{\mathrm{virt}}\right ) \\
\text{\tiny (dfn. of $q_{1\bullet}$)} & = & q_{2*}\left(
q_{1\bullet} q^*_1 \left(pr^*_1 \cv \cup pr^*_2
\prod^n_{i-1}ev^*_i (\gamma_i) \right)_{\mathbb{I}_{\cross}} \cap
pr_2^*[\mgnvbbar]^{\mathrm{virt}}\right) \\
\text{\tiny ($q_1$ is splitting)} & = & q_{2*} \left(\left(pr^*_1
\cv \cup pr^*_2 \prod
ev^*_i(\gamma_i)\right)_{\mathbb{I}_{\cross}} \cap
pr_2^*[\mgnvbbar]^{\mathrm{virt}}\right)\\
& = & st_* pr_{2*} \left(\left(pr^*_1 \cv \cup pr^*_2 \prod
ev^*_i(\gamma_i) \right)_{\mathbb{I}_{\cross}} \cap pr_2^*
[\mgnvbbar]^{\mathrm{virt}} \right)\\
\text{\tiny (dfn. of $pr_{2\bullet}$)} & = & st_* \left
(pr_{2\bullet}\left(pr^*_1 \cv \cup pr^*_2 \prod ev^*_i(\gamma_i)
\right)_{\mathbb{I}_{V}} \cap [\mgnvbbar]^{\mathrm{virt}}\right)\\
\text{\tiny (prj. fmla. for $pr_{2\bullet}$)} & = & st_* \left(
\left( pr_{2\bullet}(pr^*_1 \cv) \cup \prod ev^*_i(\gamma_i)
\right)_{\mathbb{I}_{V}} \cap [\mgnvbbar]^{\mathrm{virt}}\right)\\
\text{\tiny (Equation~\ref{eqn:flat-pb2})} & = & st_* \left (
\left((st^* p_{\bullet} \cv) \cup \prod ev^*_i(\gamma_i)
\right)_{\mathbb{I}_{V}} \cap [\mgnvbbar]^{\mathrm{virt}}\right)\\
& = & st_* \left ( (st^* p_{\bullet} \cv)_{\mathbb{I}_V} \cap
\left( (\prod ev^*_i(\gamma_i))_{\mathbb{I}_{V}} \cap
[\mgnvbbar]^{\mathrm{virt}}\right)\right)\\
\text{\tiny (dfn. of $st^*$)
}& = & st_*\left((p_{\bullet} \cv)_{st} \cap  \left( (\prod
ev^*_i(\gamma_i))_{\mathbb{I}_{V}} \cap
[\mgnvbbar]^{\mathrm{virt}}\right) \right )\\
\text{\tiny (Equation~\ref{eqn:proper-push})}& = & (p_{\bullet}
\cv)_{\mathbb{I}} \cap st_* \left( (\prod
ev^*_i(\gamma_i))_{\mathbb{I}_{V}} \cap
[\mgnvbbar]^{\mathrm{virt}}\right ).\\
\end{eqnarray*}

This completes the proof of Theorem~\ref{thm:tensor}.
\end{proof}

\section{Gromov-Witten Invariants and Tensor Products of \cft s}

Let $V$ be a smooth projective variety. The moduli space of stable
$r$-spin maps $\M_{g,n}^{1/r}(V)$ gives rise to a set of
correlators satisfying axioms analogous to those satisfied by
Gromov-Witten invariants. This follows from the fact that the
\cft\ associated to $\M_{g,n}^{1/r}(V)$ is the  tensor product of
the
Gromov-Witten \cft\ with the $r$-spin \cft.

\subsection{Axioms of Gromov-Witten classes}

\

For each $\beta \in H_2(V,\nz)$ and a stable pair of integers
$(g,n)$, define the \emph{(cohomological) correlators of
Gromov-Witten
 theory} to be linear maps $\LambdaV_{g,n,\beta}:H^\bullet(V,\nc)\to
 H^\bullet(\M_{g,n},\nc)$ such that
 \[
 \LambdaV_{g,n,\beta}(\gamma_1,\ldots,\gamma_n) = \st_*
[( \prod_{i=1}^n \ev_i^* \gamma_i)
\cap\,[\M_{g,n}(V,\beta)]^\mathrm{virt}],
\]
where $\st:\M_{g,n}(V)\to\M_{g,n}$ is the stabilization
morphism~(\ref{eq:stabmor}), $[\M_{g,n}(V,\beta)]^\mathrm{virt}$ 
is the virtual fundamental     class~(\ref{eq:fundcl}) of 
$\M_{g,n}(V)$, and $\gamma_i \in H^\bullet(V,\nc)$.

\begin{thm}[\cite{KM1}]
\label{thm:gwaxioms} Let $B(V)\subset H_2(V,\nz)$  denote the 
semigroup of numerical equivalence classes $\beta$ such that 
$\beta \cdot L \geq 0$ for all ample divisor classes $L$ in $V$. 
Let $\eta$ be the Poincar\'e pairing on $H^\bullet(V)$ and  let
$\eta_{\mu\nu} := \eta(e_\mu,e_\nu)$ be the coefficients of its
matrix  with respect to a basis $\{ e_\mu \}$ for $H^\bullet(V)$. 
Denote  by $(\eta^{\mu \nu})$ the inverse matrix of
$(\eta_{\mu\nu})$.

The collection $\{\,\LambdaV_{g,n,\beta}\,\}$ satisfies the 
following properties, called the \emph{axioms of Gromov-Witten 
invariants}. (We are using the summation convention.) 

\begin{enumerate}
\item[1.] (Effectivity) $\LambdaV_{g,n,\beta} = 0$ if $\beta\notin B(V)$.

\item[2.]
      ($S_n$-Equivariance) Each map $\LambdaV_{g,n,\beta}$ is $S_n$-equivariant,
      where $S_n$ is the symmetric group on $n$ letters.

\item[3.] (Degeneration Axioms)
\begin{enumerate}
\item
Let
\begin{equation*}
\rho_{\Gamma_\mathrm{tree}}: \M_{k,j+1} \times \M_{g-k,n-j+1} \rTo
\M_{g,n}
\end{equation*}
be the gluing map corresponding to the stable graph $$
\Gamma_\mathrm{tree}=\tcongi,$$ then the forms 
$\LambdaV_{g,n,\beta}$ satisfy the composition property 
\begin{eqnarray*}
\rho_{\Gamma_\mathrm{tree}}^*
\LambdaV_{g,n,\beta}(\gamma_1,\gamma_2\,\ldots,\gamma_n)&=&
\\
\sum_{\beta_1+\beta_2 =
\beta}\LambdaV_{k,j+1,\beta_1}(\gamma_{i_1},\ldots,\gamma_{i_j},e_\mu)
\eta^{\mu\nu} &\otimes&
\LambdaV_{g-k,n-j+1,\beta_2}(e_\nu\,\gamma_{i_{j+1}},\ldots,\gamma_{i_n})
\nonumber
\end{eqnarray*}
for all $\gamma_i \in \ch$.

\item
Let
\begin{equation*}
\rho_{\Gamma_\mathrm{loop}}: \M_{g-1,n+2} \rTo \M_{g,n}
\end{equation*}
be the gluing map corresponding to the stable graph
\begin{equation*}
\Gamma_\mathrm{loop}=\ocongi,
\end{equation*}
then
\begin{equation*}
\rho_{\Gamma_\mathrm{loop}}^*\,\LambdaV_{g,n,\beta}(\gamma_1,\gamma_2,\ldots,
\gamma_n)\,=\,
\LambdaV_{g-1,n+2,\beta}\,(\gamma_1,\gamma_2,\ldots,\gamma_n,
e_\mu, e_\nu)\,\eta^{\mu\nu}.
\end{equation*}
\end{enumerate}

\item[4.] (Identity Axiom) Let $1$  be the unit in $H^\bullet(V)$ then
\[
\LambdaV_{g,n+1,\beta}(\gamma_1,\ldots,\gamma_n,1) 
 = \pi^* \LambdaV_{g,n,\beta}(\gamma_1,\ldots,\gamma_n)
 \]
for all $\gamma_i \in H^\bullet(V)$, where
$\pi:\M_{g,n+1}\to\M_{g,n}$ is the forgetful morphism.

\item[5.] (Dimension Axiom)  Let $K_V$ denote
the canonical class of $V$. The linear map $\LambdaV_{g,n,\beta}$ 
has the grading 
\[
\left| \LambdaV_{g,n,\beta} \right| =  2 \int_\beta K_V + 2
(g-2)\dim_\nc V.
\]

\item[6.] (Divisor Axiom) Let $\alpha$ belong to $H^2(V)$, then
\[
\pi_* \LambdaV_{g,n+1,\beta}(\gamma_1,\ldots,\gamma_n,\alpha) =
 \LambdaV_{g,n,\beta}(\gamma_1,\ldots,\gamma_n) \,\int_\beta \alpha
\]
for all $\gamma_i \in H^\bullet(V)$, where
$\pi:\M_{g,n+1}(V)\to\M_{g,n}$ is the forgetful morphism.

\item[7.] (Mapping to a point)
For all $\gamma_i \in H^\bullet(V)$ we have
\[
\LambdaV_{g,n}(\gamma_1,\ldots,\gamma_n) = p_{2 *} \left[ p_1^*
(\prod_{i=1}^n \gamma_i) \cup c_d(TV \boxtimes L) \right],
\]
where $p_1 : V\times\M_{g,n}\to V$ and
$p_2:V\times\M_{g,n}\to\M_{g,n}$ are the canonical projections,
$TV$ is the tangent bundle, $L = R^1\pi_*\co_{\cc_{g,n}}$ where
$\co_{\cc_{g,n}}$ is the structure sheaf on the universal curve
$\pi: \cc_{g,n}\to \M_{g,n}$, and $d = g \dim_\nc V$ (the rank of
$TV\otimes L$).
\end{enumerate}
\end{thm}

\medskip

These properties were first presented in~\cite{KM1} and were later 
proved by various people. The theorem follows from properties of 
the virtual fundamental class, restriction properties of the 
Gromov-Witten classes,  and the geometry of the moduli space of 
stable maps into $V$. (See~\cite{Ma} for a summary of the proof.)

\begin{df}
Let $\cR$ denote the ring consisting of formal sums of expressions
$q^\beta$ with complex coefficients, where $\beta$ belongs to 
$B(V)$, subject to the relations $q^{\beta_1+\beta_2} = 
q^{\beta_1} q^{\beta_2}$. We define $\LambdaV_{g,n} : 
H^\bullet(V)^{\otimes n}\to H^\bullet(\M_{g,n},\cR)$ as 
\[
\LambdaV_{g,n} := \sum_\beta\, q^\beta \LambdaV_{g,n,\beta}.
\]
Let $\LambdaV$ denote the collection $\{\,\LambdaV_{g,n}\,\}$ and
$1$ denote the unit in $H^\bullet(V)$.
\end{df}

Because of these properties, the Gromov-Witten invariants form an 
algebra over the modular operad $H_\bullet(\M)$ or, equivalently, 
a \cft. We refer the reader to~\cite{KM1,JKV} for the definition 
of a \cft.

\begin{crl}[\cite{KM1}] The triple $(H^\bullet(V,\cR),\etaV,\LambdaV)$ forms
a \cft\ with the flat identity $1$ (over the ground ring $\cR$).
\end{crl}

\subsection{Spin \cft}

\

The \cfts\ $(\Lambda,\ch,\eta)$ whose correlators are constructed
from classes $\Lambda_{g,n,\beta}$ satisfying properties $1$ to
$7$ in Theorem \ref{thm:gwaxioms} form a special class of \cfts.
There is, however, another potential construction of \cfts.

\begin{df}
Let $r\geq 2$ be an integer and let $(\chr,\etar)$ be the $(r-1)$
dimensional $\nc$ vector space with basis $\{ e_0,\ldots,e_{r-2}
\}$ together with a metric $$\etar_{m_1,m_2} :=
\etar(e_{m_1},e_{m_2}) = \delta_{m_1+m_2,r-2}.$$ Let $\cv$ be an 
$r$-spin virtual class on $\M_{g,n}^{1/r}$ satisfying the axioms 
from~\cite[\S4.1]{JKV}. Let $$\Lambdar_{g,n} : \chr^{\otimes n} 
\to H^\bullet(\M_{g,n})$$ be defined by 
\begin{equation}\label{eq:spincorr}
\Lambdar_{g,n}(e_{m_1},\ldots,e_{m_n}) := r^{1-g} p_*
c_{g,n}^{1/r,\bm}
\end{equation}
for all nonnegative numbers $g,n$ such that $2 g - 2 + n > 0$
where $p : \M_{g,n}^{1/r,\bm}\to\M_{g,n}$.  Finally, let
$\Lambdar$ denote the collection $\{ \Lambdar_{g,n} \}$.
\end{df}

In~\cite{JKV}, following Witten \cite{W}, we constructed this
$r$-spin virtual class in genus zero for all $r$ and for all
genera when $r=2$.

\begin{thm}[{\cite[Theorem 3.8]{JKV}}]
For each integer $r\geq 2$, the collection $(\Lambdar,\chr,\etar)$
forms a \cft\ with flat identity $e_0$ and is called an
\textsl{$r$-spin \cft.}
\end{thm}

Since the space $\M_{g,n}^{1/r}$ is associated to the $r$-spin
\cft\, and the space $\M_{g,n}(V)$  is associated to Gromov-Witten
theory, it is natural to ask if there is a natural \cft\
associated to the space $\M_{g,n}^{1/r}(V,\beta)$. The answer is
yes, and this \cft\ is closely related to the other two through
the operation of tensor product.

\subsection{Tensor products of \cft s}

\

The category of cohomological field theories has a canonical
tensor product operation~(see~\cite{KMK}).

\begin{df}
Let $(\ch',\eta',\Lambda')$ and $(\ch'',\eta'',\Lambda'')$ be
\cfts. Their \textsl{tensor product} is
$(\ch'\,\otimes\,\ch'',\eta'\,\otimes\,\eta'',\Lambda)$, where
\[
\Lambda_{g,n}(v'_1\otimes v''_1,\ldots, v'_n\otimes v''_n)\,:=\,
(-1)^\sigma\,\Lambda_{g,n}'(v'_1,\ldots,
v'_n)\,\cup\,\Lambda''_{g,n} (v''_1,\ldots,v''_n)
\]
for all $v'$ in $\ch'$, $v''$ in $\ch''$, and $(-1)^\sigma$
denotes the usual sign associated to the permutation
\[ (v'_1\otimes v''_1)\otimes\cdots\otimes (v'_n\otimes
v''_n)\,\mapsto\,v'_1\otimes\ldots\otimes v'_n\otimes
v''_1\otimes\ldots\otimes v''_n.
\]
\end{df}

This reflects the fact that the diagonal map
$\M_{g,n}\,\to\,\M_{g,n}\,\times\, \M_{g,n}$ is a coproduct with
respect to the composition maps of the modular operad
$\{\,H_\bullet(\M_{g,n})\,\}$.

In the case of Gromov-Witten invariants, Behrend~\cite{B} proved
that the \cft\ arising from  $\M_{g,n}(V'\,\times\,V'')$ is the
tensor product of that arising  from $\M_{g,n}(V')$ and
$\M_{g,n}(V'')$. When restricting to genus zero, one can view this
result as a deformation of the K\"{u}nneth theorem. Similarly, it was 
shown in~\cite{JKV2} that the tensor product of an $r$-spin \cft\ 
and an $r'$-spin \cft\ can be geometrically realized  by means of 
the moduli space of $(r,r')$-spin curves. To complete this 
picture, what is missing is a description of the tensor product of 
the Gromov-Witten theory with the $r$-spin \cft.

\begin{df} \label{df:LambdaVrDef}
 Let $(H^\bullet(V,\nc),\eta_P)$ denote the cohomology of $V$ together with its
Poincar\'e pairing $\eta_P$. Let $(\chVr,\eta)$ denote the tensor
product of $(H^\bullet(V),\eta_P)$ with $(\chr,\etar)$. For each 
stable pair $(g,n)$ and $\beta \in H_2(V,\nz)$, define the 
\textsl{(cohomological) correlators} to be linear maps $$ 
\LambdaVr_{g,n,\beta}:\chVr\to H^\bullet(\M_{g,n},\nc)$$ given by 
\begin{equation}\label{eq:LambdaVrDef}
\LambdaVr_{g,n,\beta}(\gamma_1\otimes
e_{m_1},\ldots,\gamma_n\otimes e_{m_n}) = Q_* [(
\tilde{c}_{g,n}^{1/r,\bm} \prod_{i=1}^n \ev_i^* \gamma_i)
\cap\,[\M^{1/r}_{g,n}(V,\beta)]^\mathrm{virt}],
\end{equation}
where $Q:\M_{g,n}^{1/r}(V)\to\M_{g,n}$ is the morphism that
forgets both  the stable map  and the $r$-spin structure, 
$[\M^{1/r}_{g,n}(V,\beta)]^\mathrm{virt}$ is the virtual 
fundamental class of $\M^{1/r}_{g,n}(V)$, and $\gamma_i\otimes 
e_{m_i} \in \chVr$. 
\end{df}

The following theorem holds.

\begin{thm}
\label{thm:tensorproduct} Let $\LambdaVr_{g,n} : \chVr^{\otimes
n}\to H^\bullet(\M_{g,n},\cR)$, where
\[
\LambdaVr_{g,n} := \sum_\beta q^\beta \LambdaVr_{g,n,\beta}.
\]
Let $\LambdaVr$ denote the collection $\{\,\LambdaVr_{g,n}\,\}$.
The collection $(H^\bullet(V,\cR),\eta,\Lambda)$ forms a \cft\
(over the ground ring $\cR$) with flat identity $1\otimes e_0$ and
is the tensor product of the \cfts\
$(\LambdaV,H^\bullet(V,\cR),\etaV)$ and $(\Lambdar,\chr,\etar)$.
\end{thm}
\begin{proof}
This is an immediate consequence of Theorem~\ref{thm:tensor}.
\end{proof}

\medskip

The $r$-spin \cfts\ behave as though the elements of $\chr$ were
cohomology classes of fractional dimension, similar to the
orbifold cohomology classes of Chen and Ruan~\cite{CR}. There is
also no analog of the elements $\beta$ in $B(V)$ appearing in
Gromov-Witten theory. However, the theory associated to $r$-spin
maps into $V$ does satisfy analogous axioms. In particular, this
theory, like the Gromov-Witten theory, is of
qc-type~\cite{Ma,Ma2}.

\subsection{Spin Gromov-Witten invariants}

\

The classes $\LambdaVr_{g,n,\beta}$ have properties analogous to
those of Gromov-Witten invariants.

\begin{thm}
\label{thm:newaxioms} Let $(g,n)$ be a stable pair of integers. 
The collection $\{\,\LambdaVr_{g,n,\beta}\,\}$ satisfies the 
following properties: 
\begin{enumerate}
\item[1.] (Effectivity) $\LambdaVr_{g,n,\beta} = 0$ if $\beta\notin B(V)$.

\item[2.] ($S_n$-Equivariance) Each map $\LambdaVr_{g,n,\beta}$ is $S_n$-equivariant.

\item[3.] (Degeneration Axioms)
Given a basis $\{ \be_\mu \}$ for $\chVr$, let $\etaVr_{\mu\nu} :=
\etaVr(\be_\mu,\be_\nu)$
          and $(\etaVr^{\mu \nu})$ denote the inverse matrix.
\begin{enumerate}
\item
Let
\begin{equation*}
\rho_{\Gamma_\mathrm{tree}}: \M_{k,j+1} \times \M_{g-k,n-j+1} \rTo
\M_{g,n}
\end{equation*}
be defined as in Theorem~\ref{thm:gwaxioms}, then the forms 
$\LambdaVr_{g,n\beta}$ satisfy the composition property: 
\begin{eqnarray*}
\rho_{\Gamma_\mathrm{tree}}^*
\LambdaVr_{g,n,\beta}(\gamma_1,\gamma_2\,\ldots,\gamma_n)&=&
\\
\sum_{\beta_1+\beta_2 =
\beta}\LambdaVr_{k,j+1,\beta_1}(\gamma_{i_1},\ldots,\gamma_{i_j},\be_\mu)
\etaVr^{\mu\nu} &\otimes&
\LambdaVr_{g-k,n-j+1,\beta_2}(\be_\nu,\gamma_{i_{j+1}},\ldots,\gamma_{i_n})
\nonumber
\end{eqnarray*}
for all $\gamma_i \in \chVr$.
\item
Let
\begin{equation*}
\rho_{\Gamma_\mathrm{loop}}: \M_{g-1,n+2} \rTo \M_{g,n}
\end{equation*}
be defined as in Theorem~\ref{thm:gwaxioms}, then 
\begin{equation*}
\rho_{\Gamma_\mathrm{loop}}^*\,\LambdaVr_{g,n,\beta}(\gamma_1,\gamma_2,\ldots,
\gamma_n)\,=\,
\LambdaVr_{g-1,n+2,\beta}\,(\gamma_1,\gamma_2,\ldots,\gamma_n,
\be_\mu, \be_\nu)\,\etaVr^{\mu\nu}
\end{equation*}
for all $\gamma_i \in \chVr$.

\end{enumerate}

\item[4.] (Identity Axiom) Let $\bone := 1\otimes e_0$, where $1$ is the unit in
$H^\bullet(V)$ and $e_0$ the unit of $\chr$, then
\[
\LambdaVr_{g,n+1,\beta}(\gamma_1,\ldots,\gamma_n,\bone) = \pi^*
\LambdaVr_{g,n,\beta}(\gamma_1,\ldots,\gamma_n)
\]   for all $\gamma_i \in \chVr$,
where $\pi:\M_{g,n+1}(V)\to\M_{g,n}$ is the forgetful morphism.

\item[5.] (Dimension Axiom)  Let $K_V$ denote the canonical class on $V$. The
map $\LambdaVr_{g,n,\beta}$ of $\nz$-graded modules must be
homogeneous of degree
\[
\left| \LambdaVr_{g,n,\beta} \right| = 2 \int_\beta K_V + 2
(g-2)\dim_\nc V + \frac{2}{r} (r-2) (g-1).
\]

\item[6.] (Divisor Axiom) Let $\alpha\otimes e_0$ belong to $H^2(V)\otimes\chr$, then
\[
\pi_*
\LambdaVr_{g,n+1,\beta}(\gamma_1,\ldots,\gamma_n,\alpha\otimes
e_0) =
 \LambdaVr_{g,n,\beta}(\gamma_1,\ldots,\gamma_n) \,\int_\beta \alpha,
\]
for all $\gamma_i \in \chVr$, where $\pi:\M_{g,n+1}(V)\to\M_{g,n}$
is the forgetful morphism.

\item[7.] (Mapping to a Point Axiom)
\[
\LambdaVr_{g,n}(\gamma_1\otimes e_{m_1},\ldots,\gamma_n\otimes
e_{m_n}) = p_{2 *} \left[ p_1^* (\prod_{i=1}^n \gamma_i) \cup
c_d(TV \boxtimes L) \right] \cup p_*c_{g,n}^{1/r,\bm} 
\]
for all $\gamma_i \in H^\bullet(V)$, where $p_1 : V \times 
\M_{g,n} \to V$ and $p_2: V \times \M_{g,n} \to \M_{g,n}$ are the 
canonical projections, $TV$ is the tangent bundle, $L = 
R^1\pi_*\co_{\cc_{g,n}}$ where $\co_{\cc_{g,n}}$ is the structure 
sheaf on the universal curve $\pi: \cc_{g,n}\to \M_{g,n}$, and $d 
= g \dim_\nc V$ (the rank of $TV\otimes L$). Finally, 
$p:\M_{g,n}^{1/r}\to\M_{g,n}$ is the morphism forgetting the spin 
structure and $\bm = (m_1,\ldots,m_n)$. 
\end{enumerate}
\end{thm}

\begin{proof}
All axioms follow immediately from Theorems~\ref{thm:tensor} and
\ref{thm:gwaxioms}.
\end{proof}

\subsection{Potential functions and gravitational descendants}

\

Recall the potential functions associated to $\M_{g,n}(V)$. 
\begin{df}
Consider the \emph{correlation functions}
\[
\langle \tau_{a_1}(\gamma_1)\ldots\tau_{a_n}(\gamma_n)
\rangle_{g,\beta} := \int_{[\M_{g,n}(V,\beta)]^\mathrm{virt}}
\prod_{i=1}^n (\psi_i^{a_i} \ev_i^*\gamma_i)
\]
for all integers $a_1,\ldots,a_n \geq 0$ and
$\gamma_1,\ldots,\gamma_n$ in $H^\bullet(V)$. Correlation
functions such that some of the $a_i$ are nonzero are called
\emph{gravitational descendants}.

The \emph{large phase space potential (function)} associated to
$\M_{g,n}(V)$ is $$ \PhiV(\bt) := \sum_{g\geq
0}\lambda^{2g-2}\,\PhiV_g(\bt) \in
\lambda^{-2}\cR[[\lambda^2]][[t_a^\alpha ]], $$ where 
\[
\PhiV_g(\bt) := \sum_{\beta\in B(V)}\langle \exp(\bt \cdot \btau)
\rangle_{g,\beta} q^\beta
\]
and
\[
\bt\cdot\btau := \sum_{a\geq 0} \sum_{\alpha} t_a^{\alpha}
\tau_a(\ee_{\alpha})
 \]
relative to a basis $\{\,\ee_{\alpha}\,\}$ for $H^\bullet(V)$ such
that $\ee_0$ is the identity.

The \emph{small phase space potential (function)}, $\PhiV(\bx)$
where $\bx = (x^1,\ldots, x^n)$ are coordinates on $H^\bullet(V)$
relative to the basis $\{\,\ee_{\alpha}\,\}$, is obtained from
$\PhiV(\bt)$ by setting $x^\alpha := t_0^\alpha$ and $t_a^\alpha
:= 0$ for all $a\geq 1$ and all $\alpha$.

\end{df}

There are analogous potential functions associated to
$\M_{g,n}^{1/r}(V)$.

\begin{df}
Consider the \emph{correlation functions}
\[
\langle \tau_{a_1}(\gamma_1\otimes
e_{m_1})\ldots\tau_{a_n}(\gamma_n \otimes e_{m_n})
\rangle_{g,\beta} :=
\int_{[\M_{g,n}^{1/r,\bm}(V,\beta)]^\mathrm{virt}} r^{1-g}
\cvt(\bm) \prod_{i=1}^n (\psi_i^{a_i} \ev_i^*\gamma_i)
\]
    for     integers $a_1,\ldots,a_n \geq 0$,
$\gamma_1,\ldots,\gamma_n \in H^\bullet(V)$, and $e_{m_1},\ldots,
e_{m_n} \in \chr$. Correlation functions such that some of the 
$a_i$ are nonzero are called \emph{gravitational descendants}.

The \emph{large phase space potential (function)} associated to
$\M^{1/r}_{g,n}(V)$ is $$ \PhiVr(\bu) := \sum_{g\geq
0}\lambda^{2g-2}\,\PhiVr_g(\bu) \in
\lambda^{-2}\cR[[\lambda^2]][[\chVr]], $$ where
\[
\PhiVr_g(\bu) := \sum_{\beta\in B(V)}\langle \exp(\bu \cdot \btau)
\rangle_{g,\beta} q^\beta
\]
and
\[
\bu\cdot\btau := \sum_{a\geq 0} \sum_{\alpha,m} u_a^{\alpha,m}
\tau_a(\ee_{\alpha}\otimes e_m)
 \]
relative to the basis $\{\,\ee_{\alpha}\otimes e_m \,\}$ for
$\chVr$.

The \emph{small phase space potential (function)}, $\PhiVr(\by)$
where $\by$ consists of coordinates $\{ y^{\alpha,m} \}$ on
$H^\bullet(V)$ relative to the basis $\{\,\ee_{\alpha}\otimes
e_m\,\}$, is obtained from $\PhiVr(\bu)$ by setting $y^{\alpha,m}
:= u_0^{\alpha,m}$ and $u_a^{\alpha,m} := 0$ for all $a\geq 1$ and
all $\alpha,m$.
\end{df}

\begin{thm}\label{thm:unstablephi}
The small phase space potential function $\PhiVr(\by)$ is
completely determined by the potential $\PhiV(\bx)$,    the
cohomological correlators $\{\,\LambdaV_{g,n} \}$, and 
$\{\,\Lambdar_{g,n}\,\}$. 
\end{thm}

\begin{proof}
Theorem~\ref{thm:tensor} shows that the intersection numbers
$\langle \gamma_1\otimes e_{m_1} \cdots \gamma_n\otimes e_{m_n} 
\rangle_{g,n}$ are completely determined by the classes 
$\{\,\LambdaV_{g,n} \}$ and $\{\,\Lambdar_{g,n}\,\}$ if $(g,n)$ is 
stable.

              We must still 
address the unstable cases---when $(g,n) \in \{ (0,0), \ (0,1), \
(0,2), \ (1,0)\}$. But by Proposition~\ref{pr:unstab} and 
Definition \ref{r-spin}, these are always of dimension zero.

Let $\M_{g,n}^{1/r,\bm} := \coprod_i \M_{g,n}^{1/r,\bm,(i)}$ where
$\M_{g,n}^{1/r,\bm,(i)}$ are the connected components of
$\M_{g,n}^{1/r,\bm}$, and let $\pti :
\M_{g,n}^{1/r,\bm,(i)}(V,\beta) \to \M_{g,n}(V,\beta)$ be the
morphisms forgetting the $r$-spin structure.  Furthermore, let 
$\cvtim$ be $\cvt$ restricted to $\M_{g,n}^{1/r,\bm,(i)}$ and let 
us assume that $\cvtim$ is zero dimensional.  For all $\bgamma 
\otimes\be := \gamma_1\otimes e_{m_1}\cdots \gamma_n\otimes 
e_{m_n}$ in $\chVr$ we have 
\begin{eqnarray*}
 \langle \bgamma \otimes \be \rangle_{g,\beta} & = &
r^{1-g} \int (\bevt^*\bgamma\cup\cvt)\cap
[\M_{g,n}^{1/r,\bm}(V,\beta)]^\mathrm{virt} \\ & = & \sum_i
\cvtim_{g,n} r^{1-g} \int \bevt^*\bgamma \cap
[\M_{g,n}^{1/r,\bm}(V,\beta)]^\mathrm{virt} \\ & = & \sum_i
\cvtim_{g,n} r^{1-g} \int \pti^* \bev^*\bgamma \cap
\pti^*[\M_{g,n}(V,\beta)]^\mathrm{virt}  \\ & = & \sum_i
\cvtim_{g,n} r^{1-g} \int (\bev^*\bgamma)_{\pti} \cap
\pti^*[\M_{g,n}(V,\beta)]^\mathrm{virt}  \\ & = & \sum_i
\cvtim_{g,n} r^{1-g} \int \pti^* (\bev^*\bgamma \cap
[\M_{g,n}(V,\beta)]^\mathrm{virt})  \\ & = & \sum_i \cvtim_{g,n}
r^{1-g} \deg(\pti) \int \bev^*\bgamma \cap
[\M_{g,n}(V,\beta)]^\mathrm{virt} \\ & = & \sum_i \cvtim_{g,n}
r^{1-g} \deg(\pti) \langle \bgamma\, \rangle_{g,\beta}, \\
\end{eqnarray*}
where $\deg$ denotes the (orbifold) degree of $\pti$. This
completes the proof.
\end{proof}

\subsection{Gravitational descendants}
\
In this subsection, we show that when $g=0$, our constructions on
$\M_{0,n}^{1/r}(V,\beta)$ satisfy a generalization of the
so-called descent property (introduced in \cite{JKV3}). In 
particular, this property explicates the geometric origin of the 
$\psi$ classes (at least in genus zero) in the definition of the 
Gromov-Witten invariants of $V$.

It may seem curious that $\M_{g,n}^{1/r}(V,\beta)$ is defined to
be the disjoint union of  $\M_{g,n}^{1/r,\bm}(V,\beta)$, where the
$n$-tuple of nonnegative integers $\bm = (m_1,\ldots, m_n)$ is
required to satisfy $m_i \leq r-1$ for all $i=1,\ldots,n$. The 
latter restriction, however, is reasonable because of the 
isomorphism 
\[
\M_{g,n}^{1/r,\bmt}(V,\beta)\,\rTo\,\M_{g,n}^{1/r,\bmt+r\bdelta_i}(V,\beta)
\]
from Proposition \ref{prop:higherm}, where $i=1,\ldots,n$,
$\bdelta_i$ is the $n$-tuple whose $i$-th component is 1 and the
rest are zero, and $\bmt := (\mt_1,\ldots,\mt_n)$   is any
$n$-tuple of nonnegative integers.

On the other hand, in genus zero the classes $\cv(\bmt)$ change
under this identification in the following manner.

\begin{thm}(The descent property)
Let $\bmt = (\mt_1,\ldots,\mt_n)$ be an $n$-tuple of nonnegative
integers and let $\bm = (m_1,\ldots,m_n)$ be the reduction of
$\bmt \pmod r$ (i.e.,  $\bmt \equiv \bm \pmod r$ and $0 \le
m_i\leq r-1$ for  $i=1,\ldots,n$).

Let $\cvt(\bmt)$ be the top Chern class of the vector bundle
$R^1\pi_*\ce(\bmt)^*$ on  $\M_{0,n}^{1/r,\bm}$.

The following equation is satisfied on $\M_{0,n}^{1/r,\bm}$ for
all $i=1,\ldots,n$, where $\bdelta_i$ is the $n$-tuple whose
$i$-th component is 1 and the rest are zero:
\begin{equation}
\label{eq:descent} \cvt(\bmt+r\bdelta_i) =
-\frac{\mt_i+1}{r}\,\psi_i\,\cvt(\bmt).
\end{equation}

\end{thm}
\begin{proof}
The proof is identical to the case of $\M_{g,n}^{1/r}$ in
\cite{JKV3}. It follows from the long exact sequence associated to
the short exact sequence
\[
0\rTo \ce_r(\bmt) \rTo \ce_r(\bmt+r \bdelta_i)\rTo
\sigma_i^*\ce_r(\bmt)\rTo 0
\]
and the fact that
\[
\tilde{\psi}_i := c_1(\sigma^*_i \ce_r(\bmt)) =
\frac{m_i+1}{r}\,\psi_i
\]
for all $i=1,\ldots,n$, which follows from an immediate
generalization of Proposition~2.2 from \cite{JKV}.
\end{proof}

\begin{rem}\label{rem:descent}
The descent property holds on both $\M_{0,n}^{1/r}$ and
$\M_{0,n}^{1/r}(V,\beta)$, but the $\psi$ classes on
$\M_{0,n}^{1/r}(V,\beta)$ are not pullbacks of the corresponding
$\psi$ classes on $\M_{g,n}^{1/r}$---just as in the case of stable
maps, they differ by divisors that are collapsed under the
stabilization map (see~\cite{KM2,Ma}). This illustrates the fact,
alluded to in Remarks~\ref{rem:restrict-m}
and~\ref{rem:cv-m-restrict}, that when any $\tilde{m}_i$ is larger
than $r-1$, the class  $\tilde{\st}^*\cv(\bmt)$ is not equal to
the class $\cvt(\bmt)$.
\end{rem}

The previous theorem motivates the following generalization of the
small phase space potential function in genus zero.

\begin{df}

Let the $n$-tuples  $\bmt = (\mt_1,\ldots,\mt_n)$ and $\bm$ and
the class
 $\cvt(\bmt)$ on $\M_{0,n}^{1/r,\bm}$ be the same as in the previous theorem.

       Define the \emph{correlation functions}
\[
\langle \taut_{0}(\gamma_1\otimes
e_{\mt_1})\ldots\taut_{0}(\gamma_n \otimes e_{\mt_n})
\rangle_{0,\beta} :=
\int_{[\M_{0,n}^{1/r,\bm}(V,\beta)]^\mathrm{virt}} r \cvt(\bmt)
\prod_{i=1}^n \evt_i^*\gamma_i.
\]
Consider the analog of the genus zero small phase space potential
$$ \PhitVr_0(\btT) \in \cR[[\lambda^2]][[\tT^{\alpha,\mt} ]], $$ 
where 
\[
\PhitVr_0(\btT) := \sum_{\beta\in B(V)}\langle \exp(\btT \cdot
\btaut) \rangle_{0,\beta} q^\beta
\]
and
\[
\btT\cdot\btaut := \sum_{\alpha,\mt} \tT^{\alpha,\mt}
\taut_0(\ee_{\alpha}\otimes e_\mt),
 \]
where the sum runs over all $\alpha$ and all nonnegative integers
$\mt$.
\end{df}

\begin{crl}
Let $r\geq 2$ be an integer. The potential functions
$\PhitVr_0(\btT)$ and
$\PhiVr_0(\bu)$ are equal after making the assignment: 
\[
\tT^{\alpha,(a r + m)} := \frac{(-1)^a r^a}{[r(a-1)+m+1]_r}
u_a^{\alpha,m},
\]
where $a$ and $m$ are nonnegative integers such that  $m \leq r-1$
and
\[
[r(a-i)+m+1]_r := \prod_{i=1}^a ((r(a-i) + m + 1).
\]
\end{crl}

\begin{rem}
In  \cite{PoVa}, a candidate class for $\cv$ was constructed on
$\M_{g,n}^{1/r}$, which was shown to obey some of the axioms in
\cite{JKV}. In addition, it satisfies  the descent property  on 
$\M_{g,n}^{1/r}$ for nonnegative  $n$-tuples $\bmt$. The 
construction of~\cite{PoVa} can be generalized straightforwardly 
to $\M_{g,n}(V,\beta)$ to yield a class $\cvt(\bmt)$ 
             satisfying equation~(\ref{eq:descent}) and, hence, the previous
corollary for all $g$ and $n$.
\end{rem}

\subsection{The case of $r=2$}

\

In~\cite{JKV,W}, the virtual class $\cv(\bm)$ when $r=2$ was
constructed for all genera and $n$-tuples $\bm = (m_1,\ldots,m_n)$ 
with $0 \le m_i \le 1$. It was shown that the $r=2$ case reduced 
to the Gromov-Witten invariants of a point. Here we show that                                   
$2$-spin Gromov-Witten invariants are the usual Gromov-Witten 
invariants.

\begin{thm}\label{thm:requalstwo}
For a pair of nonnegative integers  $(g,n)$ and $\beta \in
H_2(V,\nz)$ let
 $\pt:\M_{g,n}^{1/2}(V,\beta)\to\M_{g,n}(V,\beta)$ be the map
      forgetting the   spin structure. For $i=1,\ldots,n$,
let $\gamma_i\otimes e_0$ belong to $\chVr$, then
\[
2^{1-g}\pt_*\left(\widetilde{c}^{1/2}(\mathbf{0}) \prod_{i=1}^n
(\evt_i^*\gamma_i) \cap
[\M_{g,n}^{1/2}(V,\beta)]^\mathrm{virt}\right) = \prod_{i=1}^n
(\ev_i^*\gamma_i) \cap [\M_{g,n}(V,\beta)]^\mathrm{virt}.
\]
Consequently, the large phase space potential functions
$\Phi^{(V,2)}(\bu)$ and $\Phi^V(\bt)$ agree after setting
$u_a^{(\alpha,0)} = t_a^\alpha$.
\end{thm}

\begin{proof}
This was proved in the case where $V$ is a point in~\cite{JKV}.
The same proof goes through here using the definition of $\cvt$
(which is now defined in the unstable range) and the fact that $
[\M_{g,n}^{1/2}(V, \beta)]^\mathrm{virt} = \tilde{p}^*
[\M_{g,n}(V,\beta)]^\mathrm{virt}$.
\end{proof}

\subsection{Genus zero and $\beta=0$}

\

Genus zero Gromov-Witten invariants of $V$ give 
 rise to the quantum cohomology of $V$, which is a certain deformation of the
cup product on $H^\bullet(V)$. The cup product itself appears as
the $\beta=0$ part of the genus zero potential function.
Similarly, the Frobenius structure associated to
$\M_{g,n}^{1/r}(V)$ can be regarded as a deformation of the
following commutative, associative product on $\chVr$.

\begin{prop}
\label{crl:betazero} Let $V$ be a smooth projective variety and
$n\geq 3$ be an integer. Let $\gamma_1,\ldots,\gamma_n$ belong to
$H^\bullet(V)$ and $e_{0},\ldots,e_{r-2}$ be the standard basis in
$\chr$, then
\[
\langle \gamma_1\otimes e_{m_1}\cdots \gamma_n\otimes e_{m_n}
\rangle_{g,\beta=0} = \int_{\M_{0,n}^{1/r,\bm}} c^{1/r}(\bm)
\,\int_V \gamma_1 \cup \ldots\cup\gamma_n.
\]
\end{prop}
\begin{proof}
This follows from the \emph{Mapping to a Point} property.
\end{proof}

\section{An Example: The small phase space potential for
$\M_{0,n}^{1/3}(\cpone)$}  \label{sec:ex}

Throughout this section let $r=3$ and $V = \cpone$. We will now
compute its genus zero small phase potential function, denoted by
\[
 \chi(\bt) := \Phi_0^{(\cpone,3)}(\bt),
\]
where $\bt$ is a set of coordinates $t^{\alpha,m}$ associated to
the basis $\{\,\tau_{\alpha,m}:=\varepsilon_\alpha\otimes e_m\,\}$
(where $\alpha = 0,1$ and $m=0,1$) for $\ch^{(\cpone,3)}$. Here
$\varepsilon_0$ is the identity element in $H^\bullet(\cpone)$ and
$\varepsilon_1$ is the element in $H^2(\cpone)$ Poincar\'e dual to
a point. The metric in this basis is
\[
\eta_{(\alpha_1,m_1),(\alpha_2,m_2)} :=
\eta(\varepsilon_{\alpha_1}\otimes e_{m_1}
,\varepsilon_{\alpha_2}\otimes e_{m_2}) =
\delta_{\alpha_1+\alpha_2,1} \delta_{m_1+m_2,1}.
\]

The potential function can be broken into two pieces:
\[
\chi(\bt) = \chi_{\beta=0}(\bt) + \Psi(\bt),
\]
where $\chi_{\beta = 0}(\bt)$ consists of only those terms
corresponding to the moduli spaces $\M_{0,n}(\cpone,0)$; while 
$\Psi(\bt)$ contains the contributions (``instanton corrections'') 
            from $\M_{0,n}(\cpone,\beta)$ where $\beta \not= 0$.
Corollary~\ref{crl:betazero} implies that
\begin{equation}
\label{eq:Phibetazero} \chi_{\beta = 0}(\bt) = \frac{1}{2} t^{1,1}
(t^{0,0})^2 + t^{0,0} t^{0,1} t^{1,0} + \frac{1}{18} t^{1,1}
(t^{0,1})^3.
\end{equation}
Theorem~\ref{thm:newaxioms} implies that
\begin{equation}
\label{eq:Psi} \Psi(\bt) = \sum_{\beta\geq 1}\sum_{n_1,n_2\geq 0}
q^\beta \frac{ (t^{0,1})^{n_1} (t^{1,0})^{n_2} (t^{1,1})^{6 \beta
+ 2 n_1- 5}} {n_1! n_2! (6 \beta + 2 n_1 - 5)!}
\cor{\tau_{0,1}^{n_1}\tau_{1,0}^{n_2}\tau_{1,1}^{6 \beta + 2 n_1-
5}}_\beta.
\end{equation}

Furthermore, Theorem~\ref{thm:newaxioms} implies that the
potential function must satisfy the WDVV equation
\begin{eqnarray*}
\frac{\partial^3 \chi(\bt)}{\partial t^{\alpha_1,m_1} \partial
t^{\alpha_2,m_2} \partial t^{\alpha_+,{m_+}}} &&
\eta^{(\alpha_+,m_+),(\alpha_-,m_-)} \frac{\partial^3
\chi(\bt)}{\partial t^{\alpha_-,{m_-}} \partial t^{\alpha_3,m_3}
\partial t^{\alpha_4,m_4}}  = \\ \frac{\partial^3
\chi(\bt)}{\partial t^{\alpha_3,m_3} \partial t^{\alpha_2,m_2}
\partial t^{\alpha_+,{m_+}}} &&
\eta^{(\alpha_+,m_+),(\alpha_-,m_-)} \frac{\partial^3
\chi(\bt)}{\partial t^{\alpha_-,{m_-}} \partial t^{\alpha_1,m_1}
\partial t^{\alpha_4,m_4}} \\
\end{eqnarray*}
for all $m_i,\alpha_i=0,1$ and $i=1,\ldots,4$, and where the
summation convention has been used.

Setting $(\alpha_1,m_1) = (1,0)$, $(\alpha_2,m_2) = (0,1)$,  and
$(\alpha_3,m_3) = (\alpha_4,m_4) = (1,1)$ in the WDVV equation and
plugging in equation (\ref{eq:Phibetazero}), we obtain
\begin{eqnarray*}
\partial^3_{1,1} \Psi =
&-& \partial^2_{0,1} \partial_{1,0} \Psi
\partial_{1,0}\partial_{1,1}^2\Psi \\ &-& \partial_{0,1}
\partial^2_{1,0} \Psi \partial_{0,1}\partial_{1,1}^2\Psi \\ &+&
\frac{1}{3} t^{0,1}\partial^2_{0,1}\partial_{1,1}\Psi \\ &+&
\partial_{0,1}^2 \partial_{1,1} \Psi
\partial_{1,0}^2\partial_{1,1}\Psi \\ &+&
(\partial_{0,1}\partial_{1,0}\partial_{1,1}\Psi)^2, \\
\end{eqnarray*}
where we have used the shorthand notation
\[
\partial^n_{\alpha, m} = \left(\frac{\partial}{\partial t^{\alpha,m}}\right)^n.
\]
Together with the Divisor Axiom in Theorem~\ref{thm:newaxioms}, we
obtain the recursion relations for $\beta = 1$ correlators
\[
\cor{\tau_{0,1}\tau_{1,1}^3}_1 = \frac{1}{3} \cor{\tau_{1,0}^2
\tau_{1,1}}_1,
\]
and, for all $n_1\geq 2$,
\[
\cor{\tau_{0,1}^{n_1}\tau_{1,1}^{2n_1+1}}_1 = \frac{n_1}{3}
\cor{\tau_{0,1}^{n_1-1} \tau_{1,1}^{2 n_1-1}}_1.
\]
These collectively imply that for all $n_1\geq 1$,
\[
\cor{\tau_{0,1}^{n_1}\tau_{1,1}^{2n_1+1}}_1 = \frac{n_1!}{3^{n_1}}
\cor{\tau_{1,0}^2\tau_{1,1}}_1.
\]
Furthermore, the tensor product property implies that
\[
\cor{\tau_{1,0}^2\tau_{1,1}}_1 = 1.
\]
Together with the Divisor Axiom, this determines all of the $\beta
= 1$ correlators.

If $\beta \geq 2$ then we obtain the following recursion relation
for all $n_1\geq 0$:
\begin{eqnarray*}
&&\cor{\tau_{0,1}^{n_1} \tau_{1,1}^{6\beta + 2 n_1 - 5}}_\beta =
 \frac{n_1 \beta^2}{3} \cor{\tau_{0,1}^{n_1-1} \tau_{1,1}^{6 \beta
+ 2 n_1 - 7}}_\beta\\
 & + & \sum ( - \beta' \beta'' {{n_1} \choose {n_1'}} {{6 \beta +
2 n_1 - 8}\choose {6 \beta' + 2 n_1' - 1}}
\cor{\tau_{0,1}^{n_1'+2}\tau_{1,1}^{6\beta'+ 2 n_1' - 1}}_{\beta'}
 \cor{\tau_{1,0}^{n_1''}\tau_{1,1}^{6 \beta''+ 2n_1'' -
5}}_{\beta''} \\ &- & (\beta')^2 {{n_1}\choose {n_1'}} {{6 \beta +
2 n_1 - 8}\choose {6 \beta' + 2 n_1' -
3}}\cor{\tau_{0,1}^{n_1'+1}\tau_{1,1}^{6\beta'+ 2 n_1' -
3}}_{\beta'} \cor{\tau_{1,0}^{n_1''+1}\tau_{1,1}^{6 \beta''+
2n_1'' - 3}}_{\beta''}\\ &+& (\beta'')^2 {{n_1}\choose {n_1'}} {{6
\beta + 2 n_1 - 8}\choose {6 \beta' + 2 n_1' -
2}}\cor{\tau_{0,1}^{n_1'+2}\tau_{1,1}^{6\beta'+ 2 n_1' -
1}}_{\beta'} \cor{\tau_{1,0}^{n_1''}\tau_{1,1}^{6 \beta''+ 2n_1''
-5}}_{\beta''}\\ &+& \beta ' \beta'' {{n_1}\choose {n_1'}} {{6
\beta + 2 n_1 - 8}\choose {6 \beta' + 2 n_1' -
4}}\cor{\tau_{0,1}^{n_1'+1}\tau_{1,1}^{6\beta'+ 2 n_1' -
3}}_{\beta'} \cor{\tau_{1,0}^{n_1''+1}\tau_{1,1}^{6 \beta''+
2n_1'' -3}}_{\beta''} ),\\
\end{eqnarray*}
where the first summation is over $\beta',\beta''\geq 1$ such that
$\beta =
 \beta' + \beta''$, and over  $n_1',n_1'' \geq 0$ such that $n_1 = n_1' +
 n_1''$.  Furthermore, we have defined
\[
\cor{\tau_{0,1}^{-1}\tau_{1,1}^{6 \beta -  7}}_\beta := 0.
\]
Together with the Divisor Axiom, these recursion relations
completely determine all of the $n$-point correlators of the
theory where $n\geq 3$.

Finally, the $0$, $1$ and $2$ point correlators (those in the
unstable range) are determined as a special case of
Theorem~\ref{thm:unstablephi}. The only nonvanishing correlators
of these types are
\[
\cor{\tau_{1,1}}_1 = \cor{\tau_{1,0}\tau_{1,1}}_1 = 1.
\]

\bibliographystyle{amsplain}

\begin{thebibliography}{10}






\bibitem{AJ}
D.\,Abramovich, T.\,Jarvis, \emph{Moduli of twisted spin curves},
in preparation, 2000.





\bibitem{AV}
D.\,Abramovich, A.\,Vistoli, \emph{Compactifying the space of
stable maps}, \texttt{math.AG/9908167}.






\bibitem{B2}
K.\,Behrend, \emph{Gromov-Witten invariants in algebraic
geometry}, Invent.\ Math.\ \textbf{127} (1997), 601--617.





\bibitem{B}
\bysame, \emph{The product formula for Gromov-Witten invariants,}
J.\ Alg.\ Geom.\ \textbf{8}, (1999),    529--541,
\texttt{alg-geom/9710014}.





\bibitem{BMa}
K.\,Behrend, Yu.\,Manin, \emph{Stacks of stable maps and
Gromov-Witten invariants.} Duke Math.\ J.\ \textbf{85} (1996), 
1--60.







\bibitem{CR} W.\,Chen, Y.\,Ruan \emph{A new cohomology theory for
       orbifold,} \texttt{math.AG/0004129}.



\bibitem{Du}
B.\,Dubrovin, \emph{Geometry of 2D topological field theories,}
``Integrable Systems and Quantum Groups,'' Lecture Notes in Math.\
\textbf{1620}, Springer-Verlag, Berlin, 1996.





\bibitem{EGA4}
A.\,Grothendieck and J.\,Dieudonn\'{e}.
\newblock {\em \'{E}l\'{e}ments de G\'{e}om\'{e}trie Alg\'{e}brique {IV}: \'{E}tude
  Locale des Sch\'{e}mas et des Morphismes de Sch\'{e}mas}, volume~28.
\newblock Publications Math\'{e}matiques IHES, 1966.





\bibitem{Hi}
N.\,Hitchin, \emph{Frobenius manifolds,} ``Gauge Theory and
Symplectic Geometry (Montreal, 1995),'' J.\,Hurtubise e.a.\
(eds.),  NATO Adv.\ Sci.\ Inst.\ Series C \textbf{488}, Kluwer
Publ., Dordrecht,  1997, 69--112.





\bibitem{J} T.\,J.\,Jarvis, \emph{Geometry of the moduli of higher spin
curves,} Internat.\ J.\ of Math.\ \textbf{11} (2000), 637--663,
 \texttt{math.AG/9809138}.





\bibitem{J2} \bysame, \emph{Torsion-free sheaves and moduli of
    generalized spin curves,} Compositio Math.\ \textbf{110} (1998),
    291--333.





\bibitem{J3} \bysame, \emph{Picard group of the moduli of higher spin
curves,} preprint, \texttt{math.AG/9908085}.





\bibitem{J4}\bysame, \emph{Compactification of the universal Picard over the moduli of stable curves,}
 Math.\ Zeitschrift.\ \textbf{235} (2000), 123--149.




\bibitem{JKV3} T.\,Jarvis, T.\,Kimura, A.\,Vaintrob, \emph{Gravitational
descendants and the moduli space of higher spin curves}, to appear
in  Advances in Algebraic Geometry Motivated by Physics,
Contemporary Mathematics, AMS, \texttt{math.AG}/0009066.




\bibitem{JKV} \bysame, \emph{Moduli spaces of higher
spin curves and integrable hierarchies}, to appear in Compositio
Math., \texttt{math.AG/9905034}.




\bibitem{JKV2} \bysame, \emph{Tensor products of Frobenius manifolds and
moduli spaces of higher spin curves}, ``Confer\'{e}nce de
Mosh\'{e} Flato 1999, Vol.\ 2,'' G. Dito and D. Sternheimer
(eds.), Kluwer (2000), 145--166, \texttt{math.AG/9911029}.




\bibitem{Ko}
M.\,Kontsevich, \emph{Intersection theory on the moduli space of
curves and the matrix Airy function,} Commun.\ Math.\ Phys.\
\textbf{147} (1992), 1--23.






\bibitem{KM1}
M.\,Kontsevich, Yu.\,I.\,Manin, \emph{Gromov-Witten classes,
quantum cohomology,
  and enumerative geometry}, Commun.\ Math.\ Phys.\ \textbf{164} (1994), 525--562.






\bibitem{KMK}
M.\,Kontsevich, Yu.\,I.\,Manin (with Appendix by R. Kaufmann),
\emph{Quantum cohomology of a product}, Invent.\ Math.\
\textbf{124} (1996), 313--340.




\bibitem{KM2}
\bysame, \emph{Relations between the correlators of the
topological sigma-model coupled to gravity,} Comm. Math. Phys.
\textbf{196} (1998), no. 2, 385--398.







\bibitem{Ma}
Yu.\,I.\,Manin, ``Frobenius manifolds, quantum cohomology, and
moduli spaces,''
 American Mathematical Society, Providence, 1999.




\bibitem{Ma2}
\bysame \emph{Three constructions of Frobenius manifolds: a
comparative study,} Asian J. Math. \textbf{3} (1999), 179--220, 
\texttt{math.AG/9801006}. 




\bibitem{Mu}
D.\, Mumford, \emph{Towards an enumerative geometry of the moduli
space of curves}, in ``Arithmetic and Geometry,'' (eds.\ M.\,Artin
and J.\,Tate), Part II, Progress in Math., Vol.\ 36,
Birkh\"{a}user, Basel (1983), 271--328.




\bibitem{PoVa}
A.~Polishchuk, A.~Vaintrob, \emph{Algebraic construction of
Witten's top Chern class},  to appear in  Advances in Algebraic
Geometry Motivated by Physics, Contemporary Mathematics, AMS, \
\texttt{math.AG/0011032}.





\bibitem{V}
A.\,Vistoli, \emph{Intersection theory on algebraic stacks and on
their moduli spaces},  Invent.\ Math. \textbf{97} (1989),
613--670.




\bibitem{W}
E.\,Witten, \emph{Algebraic geometry associated with matrix models
of two
  dimensional gravity}, Topological Methods in Modern Mathematics (Stony
Brook, NY, 1991), Publish or Perish, Houston, 1993, 235--269.





\end{thebibliography}

\providecommand{\bysame}{\leavevmode\hbox
to3em{\hrulefill}\thinspace}

\end{document}